\documentclass[10pt, a4paper, oneside,reqno]{amsart}

\makeatletter
\def\@settitle{\begin{center}%
		\baselineskip14\p@\relax
		\normalfont\LARGE\bfseries
		\@title
	\end{center}%
}

\def\section{\@startsection{section}{1}%
	\z@{.7\linespacing\@plus\linespacing}{.5\linespacing}%
	{\normalfont\large\bfseries}}

\def\subsection{\@startsection{subsection}{2}%
	\z@{.5\linespacing\@plus.7\linespacing}{.5\linespacing}%
	{\normalfont\bfseries}}

\def\@setauthors{%
  \begingroup
  \def\thanks{\protect\thanks@warning}%
  \trivlist
  \centering\footnotesize \@topsep30\p@\relax
  \advance\@topsep by -\baselineskip
  \item\relax
  \author@andify\authors
  \def\\{\protect\linebreak}%
  \authors%
  \ifx\@empty\contribs
  \else
    ,\penalty-3 \space \@setcontribs
    \@closetoccontribs
  \fi
  \endtrivlist
  \endgroup
}

\makeatother

\usepackage[table, xcdraw, usenames, dvipsnames]{xcolor}
\definecolor{darkblue}{rgb}{0.0, 0.0, 0.45}
\definecolor{darkgreen}{rgb}{0.0, 0.45, 0}
\usepackage[colorlinks	= true,
raiselinks	= true,
linkcolor	= darkblue, 
citecolor	= Mahogany,
urlcolor	= darkgreen,
pdfauthor	= {},
pdftitle	= {},
pdfkeywords	= {},
pdfsubject	= {},
plainpages	= false]{hyperref}

\pdfoutput=1
\date{\today}

\usepackage{dsfont,amsfonts,amssymb,amsmath,amsthm}
\usepackage{mathtools, mathrsfs}

\usepackage{graphicx}
\usepackage[font=small,labelfont=rm]{subcaption}
\usepackage[font=small,margin=10pt]{caption}
\usepackage{wrapfig}

\usepackage{multirow}
\usepackage{booktabs}

\usepackage{algorithm,algorithmic}

\usepackage{enumitem}
\usepackage{tikz}
\usepackage{courier} 

\usepackage{nicefrac}

\theoremstyle{plain}
\newtheorem{Thm}{Theorem}[section]

\newtheorem{Lem}[Thm]{Lemma}

\newtheorem{Rem}[Thm]{Remark}

\DeclareMathOperator*{\argmin}{\arg\!\min}

\newcommand{\R}{\mathbb{R}}

\newcommand{\ra}{\rightarrow}

\newcommand{\ol}[1]{\overline{#1}}

\newcommand{\wt}{\widetilde}
\newcommand{\wh}{\widehat}
\newcommand{\Let}{\coloneqq}

\newcommand{\tr}{^{\top}}

\newcommand{\norm}[1]{\left\Vert #1 \right\Vert}
\newcommand{\inner}[2]{\left\langle #1, #2 \right\rangle }

\newcommand{\opt}{^\star}

\newcommand{\EE}{\mathds{E}}

\newcommand{\set}[1]{\mathcal{#1}}

\DeclareMathOperator{\ord}{\mathcal{O}}

\graphicspath{{Fig/}}

\usepackage[square,numbers]{natbib}
\usepackage{bm}

\newcommand{\bo}{\mathbf{T}}
\newcommand{\sbo}{\wh{\mathbf{T}}}
\newcommand{\e}{\bm{1}}
\newcommand{\z}{\bm{0}}

\title[Rank-One Modified Value Iteration]{Rank-One Modified Value Iteration}

\author[]{Arman~S.~Kolarijani$^1$, Tolga~Ok$^1$, Peyman~Mohajerin~Esfahani$^{1,2}$, and Mohamad Amin Sharif~Kolarijani$^1$}
\thanks{$^1$Delft Center for Systems and Control, 
Delft University of Technology, The Netherlands. 
$^2$Department of Mechanical and Industrial Engineering, University of Toronto, Canada. 
Correspondence to: Arman S. Kolarijani $<$\texttt{a.sharifikolarijani@tudelft.nl}$>$.}
\thanks{This work was partially supported by the Horizon Europe Pathfinder Open project RELIEVE-101099481 and by the European Research Council (ERC) project TRUST-949796.} 
\thanks{The authors would like to thank the reviewers of ICML 2025 for their useful comments.}

\begin{document}

\begin{abstract}
In this paper, we provide a novel algorithm for solving planning and learning problems of Markov decision processes. 
The proposed algorithm follows a policy iteration-type update by using a rank-one approximation of the transition probability matrix in the policy evaluation step. 
This rank-one approximation is closely related to the stationary distribution of the corresponding transition probability matrix, 
which is approximated using the power method. 
We provide theoretical guarantees for the convergence of the proposed algorithm to optimal (action-)value function with the same rate and computational complexity as the value iteration algorithm in the planning problem and as the Q-learning algorithm in the learning problem. 
Through our extensive numerical simulations, however, we show that the proposed algorithm consistently outperforms first-order algorithms and their accelerated versions for both planning and learning problems.


\end{abstract}

\maketitle

\section{Introduction}
\label{sec:intro}



Value iteration (VI) and policy iteration (PI) lie at the heart of most if not all algorithms for optimal control of Markov decision processes (MDPs) in both cases of the planning problem (i.e., with access to the true model of the MDP) and the reinforcement learning problem (i.e., with access to samples of the MDP)~\cite{sutton2018reinforcement,bertsekas2023course}. 
Their widespread application stems from their simple implementation and straightforward combination with function approximation schemes such as neural networks. 
Both VI and PI are iterative algorithms that ultimately find the fixed-point of the Bellman (optimality) operator~$\bo$. 
To be precise, for $\gamma$-discounted, finite state-action MDPs, the value function $\bm{v}_k$ at iteration~$k$ is given by
\begin{equation}\label{eq:VI-PI}
\bm{v}_{k+1} = \bm{v}_k + \bm{G}_k \big( \bo (\bm{v}_k) -  \bm{v}_k\big) , \quad k=0,1,\ldots,    
\end{equation}
with $\bm{G}_k = \bm{I}$ in the VI algorithm and $\bm{G}_k = (\bm{I} - \gamma \bm{P}_k)^{-1}$ in the PI algorithm, where $\bm{I}$ is the identity matrix and $\bm{P}_k$ is the state transition probability matrix of the MDP under the greedy policy with respect to~$\bm{v}_k$. 
Both algorithms are guaranteed to converge to the optimal value function and control policy; see, e.g.,~\cite[Thms.~6.3.3 and 6.4.2]{puterman2014markov}. 
When it comes to computational complexity, one observes a trade-off between the two algorithms: VI has a lower per-iteration complexity compared to PI, while PI converges in a fewer number of iterations compared to VI. 
The faster convergence of PI is partially explained by its second-order nature which leads to a local quadratic convergence rate~\cite{puterman1979convergence, bertsekas2022lessons, gargiani2022dynamic}, compared to the linear convergence rate of VI~\cite[Thms.~6.3.3]{puterman2014markov}. 

This trade-off between the two algorithms has been the motivation for much research aimed at improving the convergence rate of VI and/or the per-iteration complexity of PI. 
One of the first improvements is the \emph{Relaxed VI} algorithm~\cite{kushner1971accelerated,porteus1978accelerated} which allows for a greater step size compared to the standard VI algorithm. 
More importantly, the correspondence between VI and gradient decent algorithm and between PI and Newton method~\cite{grand2021convex, kolarijani2023optimization} has led to a large body of research adapting ideas such as accelerated methods and quasi-Newton methods from the optimization literature for developing modified versions of VI and PI. 
For instance, the combination of the VI algorithm with Nesterov acceleration~\cite{nesterov1983method} and Anderson acceleration~\cite{anderson1965iterative} have been explored in~\cite{goyal2019first} and~\cite{zhang2020globally}, respectively, for solving the planning problem. 
More recently, Halpern's anchoring acceleration scheme~\cite{halpern1967fixed} has been used to introduce the \emph{Anchored VI} algorithm \cite{lee2024accelerating} which in particular exhibits a $\ord(1/k)$-rate for large values of discount factor and even $\gamma = 1$. 
In the case of the learning problem, \emph{Speedy Q-Learning}~\cite{ghavamzadeh2011speedy}, \emph{Momentum Q-Learning}~\cite{weng2020momentum}, and \emph{Nesterov Stochastic Approximation}~\cite{devraj2019matrix} are among the algorithms that use the idea of momentum to achieve a better rate of convergence compared to standard Q-learning (QL). 
The \emph{Quasi-Policy Iteration/Learning} algorithms~\cite{kolarijani2023optimization} are, on the other hand, an example of using the idea of quasi-Newton methods for developing algorithms for optimal control of MDPs. 
Another class of modified VI algorithms is the \emph{Generalized Second-Order VI} algorithm~\cite{kamanchi2021generalized} which applies the Newton method on a \emph{smoothed} version of the Bellman operator. 
Tools and techniques from linear algebra have also been exploited to modify the VI algorithm, particularly for policy evaluation. 
The \emph{Operator Splitting VI} algorithm \cite{rakhsha2022operator} is an example that exploits the matrix splitting method for solving the linear equation corresponding to policy evaluation for a given ``cheap-to-access'' model of the underlying MDP. 
Recently, in~\cite{lee2024deflated}, the authors combined the matrix splitting method with the matrix deflation techniques for removing the dominant eigenstructure of the transition probability matrix to speed up the policy evaluation.

\noindent\textbf{Contribution.} 
In this paper, we propose a novel algorithm that modifies the VI algorithm by incorporating a computationally efficient PI-type update rule.
To this end, we consider the update rule~\eqref{eq:VI-PI} with a matrix gain of the form $\bm{G}_k = (\bm{I} - \wt{\bm{P}}_k)^{-1}$, where $\wt{\bm{P}}_k$ is a \emph{rank-one approximation} of $\bm{P}_k$. 
To be precise, we consider the approximation~$\wt{\bm{P}}_k = \e \bm{d}_k\tr$, where $\bm{d}_k = \bm{P}_k\tr \bm{d}_k$ is a \emph{stationary} distribution of $\bm{P}_k$ and $\e$ is the all-one vector. 
The proposed algorithm then uses the \emph{power method} for approximating $\bm{d}_k$ iteratively using the true matrix~$\bm{P}_k$ in the planning problem and its sampled version in the learning problem. 
In particular, 
\begin{itemize}
    \item[(1)] we propose the \emph{Rank-one VI} (\emph{R1-VI}) Algorithm~\ref{alg:R1-VI} as the modified VI algorithm for solving the planning problem and prove its convergence to the optimal value function (Theorem~\ref{thm:R1VI-conv}); 
    \item[(2)] we propose the \emph{Rank-one QL} (\emph{R1-QL}) Algorithm~\ref{alg:R1-QL} as the modified QL algorithm for solving the learning problem and prove its convergence  to the optimal Q-function (Theorem~\ref{thm:R1Ql-conv});
    \item[(3)] we compare the proposed R1-VI and R1-QL algorithms with the state-of-the-art algorithms for solving planning and learning problems of MDPs and show the empirically faster convergence of the proposed algorithms compared to the ones with the same per-iteration computational complexity (i.e., the first-order algorithms and their accelerated versions).
\end{itemize}

\noindent\textbf{Paper organization.} 
In Section~\ref{sec:control problem}, we provide the necessary background and the problem definition along with the standard VI and QL algorithms for solving the planning and learning problems of MDPs. 
The proposed R1-VI algorithm for solving the planning problem and its analysis are discussed in Section~\ref{sec:R1-VI}. 
Section~\ref{sec:R1QL} presents the R1-QL algorithm for solving the learning problem and its analysis. 
In Section~\ref{sec:num}, we provide the results of our extensive numerical simulations and compare the proposed algorithms with a range of existing algorithms for solving the optimal control problem of MDPs. 
Finally, some limitations of the proposed algorithms and future research directions are discussed in Section~\ref{sec:lim}. 
All the technical proofs are provided in Appendix~\ref{app:proofs}. 

\noindent\textbf{Notations.} The set of real numbers is denoted by $\R$. For a vector~$\bm{v} \in \R^n$, we use~$\bm{v}(i)$ and~$[\bm{v}](i)$ to denote its~$i$-th element. 
Similarly,~$\bm{M}(i,j)$ and~$[\bm{M}](i,j)$ denote the element in~$i$-th row and~$j$-th column of the matrix~$\bm{M} \in \R^{m \times n}$. 
We use $\inner{\bm{v}}{\bm{u}} = \sum_{i=1}^n \bm{v}(i) \cdot \bm{u}(i)$ to denote the the inner product of the two vectors $\bm{v},\bm{u} \in \R^n$. 
$\norm{\bm{v}}_1 = \sum_{i=1}^n |\bm{v}(i)|$, $\norm{\bm{v}}_2 = \sqrt{\inner{\bm{v}}{\bm{v}}}$, and $\norm{\bm{v}}_{\infty} = \max_{i=1}^n |\bm{v}(i)|$ denote the 1-norm, 2-norm, and~$\infty$-norm of the vector~$\bm{v} \in \R^n$, respectively. 
We use $\rho(\bm{M})$ to denote the spectral radius (i.e., the largest eigenvalue in absolute value) of a square matrix $\bm{M}\in \R^{n \times n}$.
Given a set $\set{X}$, $\Delta(\set{X})$ denotes the set of probability distributions on $\set{X}$. 
Let~$x\sim P$ be a random variable with distribution~$P \in \Delta(\set{X})$. 
We use~$\wh{x} \sim P$ to denote \emph{a sample of the random variable~$x$} drawn from the sample space~$\set{X}$ of~$x$ according to the distribution~$P$. 
We use $\e$, $\z$, and $\bm{I}$ to denote the all-one vector, the all-zero vector, and the identity matrix, respectively, with their dimension being clear from the context.

\section{Optimal control of MDPs}\label{sec:control problem} 

Consider a finite MDP~$(\set{S},\set{A},P,c,\gamma)$. Here,~$\set{S} := \{1,2,\ldots,n \}$ and~$\set{A} := \{1,2,\ldots,m \}$ are the \emph{state} and \emph{action spaces}, respectively. 
The \emph{transition kernel}~$P: \set{S}\times\set{A}\ra \Delta(\set{S})$ is the conditional probability~$ P(s^+|s,a)$ of the transition to state $s^+$ given the current state-action pair $(s,a)$. 
The function~$\bm{c} \in \R^{|\set{S}\times\set{A}|}=\R^{nm} $, bounded from below, represents the \emph{stage cost}~$\bm{c}(s,a)$ of taking the control action~$a$ while the system is in state~$s$. 
And, $\gamma \in (0,1)$ is the \emph{discount factor} which can be seen as a trade-off parameter between short- and long-term costs. 

A \emph{control policy}~$\pi:\set{S}\ra\set{A}$ is a mapping from states to actions. Fix policy~$\pi$. 
For the corresponding Markov chain under the policy~$\pi$ we define:
\begin{itemize}
    \item[(i)] the \emph{state transition probability matrix}~$\bm{P}^\pi \in \R^{|\set{S}|\times|\set{S}|} = \R^{n\times n}$ where~$\bm{P}^\pi(s,s^+) = P\big(s^+ | s, \pi(s)\big)$ for $s,s^+ \in \set{S}$;
    \item[(ii)] the \emph{state-action transition probability matrix}~$\ol{\bm{P}}^\pi \in \R^{|\set{S}\times\set{A}|\times|\set{S}\times\set{A}|} = \R^{(nm)\times(nm)}$ where~
    \begin{equation*}
        \ol{\bm{P}}^\pi\big((s,a),(s^+,a^+)\big) =\left\{ \begin{array}{ll}
            P(s^+ | s, a) & \text{if } a^+ = \pi(s^+),  \\
            0 & \text{otherwise},
        \end{array}\right.
        \quad  \text{for } (s,a),(s^+,a^+) \in \set{S}\times\set{A}.
    \end{equation*}
    \item[(iii)] the stage cost~$ \bm{c}^{\pi} \in \R^{|\set{S}|} = \R^n$ where $\bm{c}^{\pi}(s) = \bm{c}\big(s,\pi(s)\big)$ for $s\in \set{S}$.
\end{itemize}

Under policy~$\pi$, the \emph{value function}~$\bm{v}^{\pi} \in \R^{|\set{S}|} = \R^n$ is the \emph{expected discounted cost} endured by following policy~$\pi$ over an infinite-horizon trajectory, that is,
\begin{align*}
    \bm{v}^{\pi}(s) &\Let \mathbb{E}_{s_{t+1}\sim \bm{P}^{\pi}(s_t,\cdot)}\left[\sum_{t=0}^\infty\gamma^t \bm{c}^{\pi}(s_t)\,\middle|\, s_0 = s\right], \quad \forall s \in \set{S}.
\end{align*}
The \emph{action-value function}~$\bm{q}^{\pi} \in \R^{|\set{S} \times \set{A}|} = \R^{nm}$ (or the so-called~\emph{Q-function}) under policy~$\pi$ is defined as 
\begin{equation*}
    \bm{q}^{\pi}(s,a) \Let \bm{c}(s,a) + \gamma\  \mathbb{E}_{s^+\sim P(\cdot|s,a)}\left[ \bm{v}^{\pi}(s^+) \right], \quad \forall (s,a) \in\set{S}\times\set{A}.
\end{equation*}
These functions can be shown to satisfy the fixed-point equations~\citep[Thm. 6.1.1]{puterman2014markov}
\begin{equation}\label{eq:pol_eval}
\bm{v}^{\pi} = \bm{c}^{\pi} + \gamma  \bm{P}^{\pi} \bm{v}^{\pi},\quad \bm{q}^{\pi} = \bm{c} + \gamma  \ol{\bm{P}}^{\pi} \bm{q}^{\pi}.  
\end{equation}

The problem of interest is to \emph{control} the MDP in a manner that the expected, discounted, infinite-horizon cost is minimized. To do so, one aims to find the \emph{optimal policy}~$\pi^*$ such that for any policy~$\pi$,
\begin{equation*}
    \bm{v}\opt(s) \Let  \bm{v}^{\pi\opt} (s) \leq \bm{v}^{\pi} (s), \quad \forall s \in\set{S},
\end{equation*}
or, equivalently, 
\begin{equation*}
    \bm{q}\opt(s,a) \Let  \bm{q}^{\pi\opt} (s,a) \leq \bm{q}^{\pi} (s,a), \quad \forall (s,a) \in\set{S}\times\set{A}.
\end{equation*} 
The optimal value function can be characterized as the solution to the fixed-point equation $\bm{v}\opt  = \bo (\bm{v}\opt)$, where $\bo:\R^{|\set{S}|}\ra\R^{|\set{S}|}$ is the  so-called \emph{Bellman} (\emph{optimality}) \emph{operator} defined as follows:
\begin{equation}\label{eq:Bellman_operator}
    [\bo (\bm{v}\opt)](s) \Let \min_{a \in\set{A}}\left\{ \bm{c}(s,a) + \gamma \ \EE_{s^+\sim P(\cdot|s,a)} \left[  \bm{v}\opt(s^+)\right] \right\}, \quad \forall s \in \set{S}.
\end{equation}
The operator~$\bo$ is a~$\gamma$-contraction in the~$\infty$-norm (i.e.,~$\norm{\bo(\bm{v}) - \bo(\bm{w})}_{\infty} \leq \gamma \norm{\bm{v}-\bm{w}}_{\infty}$ for all~$\bm{v},\bm{w}\in\R^n$)~\citep[Prop.~6.2.4]{puterman2014markov}. This contraction property is essentially the basis for the VI algorithm, introduced in~\eqref{eq:VI}.
\begin{equation}\label{eq:VI}
    \begin{array}{l}
        \texttt{intialize } \bm{v}_0 \in \R^{|\set{S}|} = \R^n \\
        \texttt{for } k = 0,1,\ldots  \\
        \hspace{0.5cm} \bm{v}_{k+1}(s) = [\bo (\bm{v}_k)](s), \quad \forall s \in\set{S}\\
        \texttt{endfor }
    \end{array}
\end{equation}
From the Banach fixed-point theorem (see, e.g., \citep[Thm.~6.2.3]{puterman2014markov}), the VI algorithm converges to~$\bm{v}\opt$ with a linear rate~$\gamma$. 
Correspondingly, one can derive the fixed-point characterization $\bm{q}\opt = \ol{\bo} (\bm{q}\opt)$ of the optimal Q-function, where~\citep[Fact~3]{szepesvari2009algorithms}
\begin{align*}
     [\ol{\bo} (\bm{q}\opt)](s,a) \Let \bm{c}(s,a) 
     + \gamma\ \EE_{s^+\sim P(\cdot|s,a)}\left[ \min_{a^+\in\set{A}}  \bm{q}\opt(s^+,a^+)  \right], \quad \forall (s,a) \in\set{S}\times\set{A}.
\end{align*}
This characterization is particularly useful when one only has access to \emph{samples}~$\wh{s}^{~+}\sim P(\cdot|s,a)$ of the next state~$s^+$ for each state-action pair $(s,a) \in \set{S}\times\set{A}$ (and not to the true transition probability kernel of the MDP). 
In particular, let us define the \emph{empirical Bellman operator}~$\sbo:\R^{|\set{S}\times\set{A}|}\times\set{S}\ra\R^{|\set{S}\times\set{A}|}$ as follows:
\begin{align*}
      [\sbo (\bm{q},s^+)](s,a) & \Let \bm{c}(s,a) + \gamma \ \min_{a^+\in\set{A}}  \bm{q}(s^+,a^+), \quad \forall (s,a,s^+) \in\set{S}\times\set{A}\times\set{S}.
\end{align*}
A classic algorithm for the learning problem is then the (\emph{synchronous}) QL algorithm~\cite{watkins1992q,kearns1998finite}, given in~\eqref{eq:QL}. 
\begin{equation}\label{eq:QL}
    \begin{array}{l}
        \texttt{intialize } \bm{q}_0 \in  \R^{|\set{S} \times \set{A}|} = \R^{nm} \\
        \texttt{for } k = 0,1,\ldots  \\
        \hspace{0.4cm} \texttt{for } (s,a) \in\set{S}\times\set{A}\\
        \hspace{.8cm} \wh{s}_k^{~+}\sim P(\cdot|s,a)\\
        \hspace{.8cm} \delta_k = \bm{q}_k(s,a) - [\sbo (\bm{q}_k,\wh{s}_k^{~+})](s,a) \\ 
        \hspace{.8cm} \bm{q}_{k+1}(s,a) = \bm{q}_k(s,a) -  \lambda_k \delta_k\\
        \hspace{0.4cm} \texttt{endfor }\\
        \texttt{endfor }
    \end{array}
\end{equation}
Above,~$\lambda_k \geq 0$ are the step-sizes. 
The QL algorithm is also guaranteed to converge to~$\bm{q}\opt$ almost surely given that the step-sizes satisfy the \emph{Robbins–Monro conditions} (i.e., $\sum_{k=0}^{\infty} \lambda_k =\infty$ and $\sum_{k=0}^{\infty} \lambda^2_k <\infty $)~\cite{tsitsiklis1994asynchronous,jaakkola1993convergence}. 
In particular, with a polynomial step-size $\lambda_k = 1/(1+k)^{\omega}$ with $\omega \in (\nicefrac{1}{2},1)$, QL outputs an $\epsilon$-accurate Q-function with high probability after 
$ \tilde{\ord} (\tau^{-4/\omega} \cdot \epsilon^{-2/\omega} +\tau^{1/(1-\omega)}) $ iterations of synchronous sampling with $\tau = 1-\gamma$~\cite{even2003learning}, while with a re-scaled linear step-size $\lambda_k = 1/(1+\tau k)$, QL has been shown to require 
$ \tilde{\ord} (\tau^{-5} \cdot \epsilon^{-2}) $ iterations of synchronous sampling for the same performance~\cite{wainwright2019stochastic}.

\section{Rank-one value iteration (R1-VI)}\label{sec:R1-VI} 

Another well-known algorithm for solving the planning problem in MDPs is the PI algorithm. 
In order to provide this algorithm in a compact form, let us first introduce the notion of \emph{greedy policy}. Given a value function~$\bm{v} \in \R^n$, the greedy policy with respect to~$\bm{v}$, denoted by~$\pi^{\bm{v}}:\set{S}\ra\set{A}$, is
\begin{equation}\label{eq:greedy_pol}
    \pi^{\bm{v}}(s) \in \argmin_{a\in\set{A}} \EE_{s^+\sim P(\cdot|s,a)}\!\! \left[ \bm{c}(s,a) + \gamma  \bm{v}(s^+)\right], \quad \forall s \in \set{S}.
\end{equation}
The PI algorithm is then summarized in~\eqref{eq:PI}. 
\begin{equation}\label{eq:PI}
    \begin{array}{l}
        \texttt{intialize } \pi_0:\set{S}\ra\set{A} \\ 
        \texttt{for } k = 0,1,\ldots  \\ 
        \hspace{0.5cm}  \bm{v}_{k} = \bm{v}^{\pi_k} \quad\;\;  \text{[policy evaluation -- eq.~\eqref{eq:pol_eval}]}  \\ 
        \hspace{0.5cm}  \pi_{k+1} =  \pi^{\bm{v}_k} \;\;  \text{[policy improvement -- eq.~\eqref{eq:greedy_pol}]} \\ 
        \texttt{endfor }
    \end{array}
\end{equation}
The algorithm to be proposed in this section is based on an alternative representation of the iterations of the PI algorithm; see also \citep[Prop.~6.5.1]{puterman2014markov}.

\begin{Lem}[Policy iteration]\label{lem:PI} 
Each iteration of the PI algorithm~\eqref{eq:PI} equivalently reads as
\begin{equation}\label{eq:PI_1}
    \bm{v}_{k+1} = \bm{v}_{k} + (\bm{I} - \gamma \bm{P}_k)^{-1} \big(\bo(\bm{v}_{k}) - \bm{v}_{k} \big),
\end{equation}
where~$\bm{P}_k \Let \bm{P}^{\pi^{\bm{v}_k}}$ is the state transition probability matrix of the MDP under the greed policy~$\pi^{\bm{v}_k}$.
\end{Lem}
The PI algorithm outputs the optimal policy in a finite number of iterations~\citep[Thm.~6.4.2]{puterman2014markov}.
Moreover, the algorithm has a local \emph{quadratic} rate of convergence when initiated in a small enough neighborhood around the optimal solution~\cite{bertsekas2022lessons, gargiani2022dynamic}. 
The faster convergence of PI compared to VI comes however with a higher per-iteration computational complexity: The per-iteration complexities of VI and PI are~$\ord(n^2m)$ and~$\ord(n^2m+n^3)$, respectively. 
The extra~$\ord(n^3)$ complexity is due to the policy evaluation step, i.e., solving a linear system of equations; see also the matrix inversion in the characterization~\eqref{eq:PI_1}. 
To address this issue, we propose to use a low-rank approximation of~$\bm{P}_k$ instead. 
Such approach allows us to approximate~$(\bm{I} - \gamma \bm{P}_k)^{-1}$ with a reduced computational cost by using the Woodbury formula~\cite{hager1989updating}. 
To be precise, we propose the \emph{rank-one VI} (\emph{R1-VI}) \emph{algorithm}
\begin{equation}\label{eq:R1-VI_first}
    \bm{v}_{k+1} = \bm{v}_{k} + (\bm{I} - \gamma \wt{\bm{P}}_k)^{-1} \big(\bo(\bm{v}_{k}) - \bm{v}_{k} \big),
\end{equation}
where 
\begin{equation}\label{eq:R1-approx}
  \wt{\bm{P}}_k = \e \bm{d}_k\tr,   
\end{equation}
is a \emph{rank-one} approximation of the true transition probability matrix~$\bm{P}_k$ at iteration~$k$, 
with $\bm{d}_k \Let \bm{d}_{\pi_k} \in \Delta(\set{S})$ being a stationary distribution of the greedy policy~$\pi_k = \pi^{\bm{v}_k}$, i.e., a solution of $ \bm{d}_k\tr \bm{P}_k= \bm{d}_k\tr$. 
Under certain conditions, \eqref{eq:R1-approx} is indeed ``the best'' rank-1 approximation of $\bm{P}_k$: 

\begin{Lem} [Rank-1 approximation]\label{lem:unique_sol}
Assume that the transition probability matrix~$\bm{P}_k$ is ergodic (i.e., irreducible and aperiodic). Then,  
\begin{equation}\label{eq:QPI-R1 approx}
\begin{aligned}
    \e \bm{d}_k\tr = &\argmin\limits_{\bm{P} \in \R^{n\times n}} \  \rho(\bm{P} - \bm{P}_k)\\
    & \ \emph{\text{s.t.}}\ \bm{P} \geq 0,\  \bm{P}\e=\e,\ \emph{\text{rank}}(\bm{P}) = 1.
\end{aligned}
\end{equation}
where $\bm{d}_k$ is the unique stationary distribution of $\bm{P}_k$. 
That is, $\wt{\bm{P}}_k = \e \bm{d}_k\tr$ is the best rank-1 approximation of $\bm{P}_k$ in terms of the spectral radius.
\end{Lem}

Using the Woodbury formula, we then have 
\[
    (\bm{I} - \gamma \wt{\bm{P}}_k)^{-1} =  (\bm{I} - \gamma \e \bm{d}_k\tr)^{-1} = \bm{I} + \tfrac{\gamma}{1-\gamma}\e \bm{d}_k\tr,
\]
and hence the R1-VI update~\eqref{eq:R1-VI_first} reads as
\begin{equation}\label{eq:R1-VI_gen}
    \bm{v}_{k+1} = \bo(\bm{v}_k) + \tfrac{\gamma}{1-\gamma} \big\langle \bm{d}_k, \bo(\bm{v}_k) - \bm{v}_k \big\rangle \e.
\end{equation}

Next to be addressed is the computation of the vector~$\bm{d}_k$. 
Considering the fact that $\bm{d}_k$ is a left eigenvector of $\bm{P}_k$ corresponding to the eigenvalue $1$, we can use the power method~\citep[Sec.~7.3]{golub2013matrix} to compute it as follows
\begin{equation}\label{eq:PowerMethod}
\bm{f} = \bm{P}_k\tr \bm{d}_k^{(i)},\ \bm{d}_k^{(i+1)} = \frac{\bm{f}}{\norm{\bm{f}}_1}  , \quad i =0,1,\ldots, I-1,    
\end{equation}
with some initialization $\bm{d}_{k}^{(0)} \in \Delta(\set{S})$ and $I \in \{1,2,\ldots \}$. We then use $\bm{d}_{k} = \bm{d}_{k}^{(I)}$ in the update rule~\eqref{eq:R1-VI_gen}. 
We note that the normalization is only to avoid the accumulation of numerical errors; to see this note that for the row stochastic matrix $\bm{P}_k$, we have $\bm{P}_k\tr \bm{d} \in \Delta(\set{S})$ for any $\bm{d} \in \Delta(\set{S})$. 
Under the assumptions of Lemma~\ref{lem:unique_sol}, the preceding iteration converges linearly to the unique stationary distribution with a rate equal to the second largest eigenvalue modulus of $\bm{P}_k$ ~\citep[Thm.~3.4.1]{gallager2011discrete}. 

\begin{algorithm}[t]
\begin{small}
   \caption{Rank-One Value Iteration (R1-VI)} 
   \label{alg:R1-VI}
\begin{algorithmic}[1]
	\REQUIRE transition kernel~$P: \set{S}\times\set{A}\ra \Delta(\set{S})$; cost function~$\bm{c} \in \R^{|\set{S}\times\set{A}|}$; discount factor~$\gamma \in (0,1)$;  
	\ENSURE optimal value function~$v\opt$
  	
	\vspace{.2cm}
  	
    \STATE initialize: $\bm{v}_0\in  \R^{|\set{S}| },\; \bm{d}_{-1}\in \Delta(\set{S})$;
    
    \vspace{.1cm}
    
    \FOR{$k = 0, 1, 2, \ldots$} 
        \STATE $\bm{P}_k = \z \in \R^{|\set{S}| \times |\set{S}|}$;
        \FOR{$s\in\set{S}$}
            \vspace{.1cm}
            \STATE $\left\{\begin{array}{l}
                a_k \in \argmin\limits_{a \in\set{A}}\left\{ \bm{c}(s,a) + \gamma \ \EE_{s^+} \left[  \bm{v}_k(s^+)\right]\right\}, \\
                {[\bo (\bm{v}_k)]}(s) = \min\limits_{a \in\set{A}}\left\{ \bm{c}(s,a) + \gamma \ \EE_{s^+} \left[  \bm{v}_k(s^+)\right]\right\};
            \end{array}\right.$
            \vspace{.1cm}
            \STATE $\bm{P}_k(s,s^+) = P(s^+|s,a_k),\ \forall s^+ \in \set{S}$;
        \ENDFOR

        \STATE $\bm{f} = \bm{P}_k\tr \bm{d}_{k-1}$; $\bm{d}_k = \bm{f} / \norm{\bm{f}}_1$; 
        \STATE $\bm{v}_{k+1} = \bo(\bm{v}_k) + \frac{\gamma}{1-\gamma} \big\langle \bm{d}_k, \bo(\bm{v}_k) - \bm{v}_k \big\rangle \e$; 
           
    \ENDFOR
   
\end{algorithmic}
\end{small}
\end{algorithm}

The complete description of the proposed R1-VI algorithm is provided in Algorithm~\ref{alg:R1-VI}. 
We note that Algorithm~\ref{alg:R1-VI} includes a single iteration (i.e., $I=1$) of the power method in~\eqref{eq:PowerMethod} initialized by~$ \bm{d}_{k}^{(0)} = \bm{d}_{k-1}$. 
The reason for this choice is that the greedy policy~$\pi^{\bm{v}_k}$ and hence the corresponding transition matrix $\bm{P}_k$ usually stays the same over multiple iterations $k$ of the algorithm in the value space. 
This means that the algorithm effectively performs multiple iterations of the power method. 
Despite using this \emph{approximation} $\bm{d}_k$ of the stationary distribution of $\bm{P}_k$ with a single iteration of the power method, the proposed algorithm can be shown to converge.

\begin{Thm}[Convergence of R1-VI]\label{thm:R1VI-conv}
    The iterates~$\bm{v}_k$ of the R1-VI Algorithm~\ref{alg:R1-VI} converge to the optimal value function~$\bm{v}\opt=\bo(\bm{v}\opt)$ with at least the same rate as VI, i.e., with linear rate~$\gamma$.
\end{Thm}

Let us also note that the per-iteration complexity of the proposed R1-VI Algorithm~\ref{alg:R1-VI} is $\ord(n^2m)$, i.e., the same as that of VI. 
We finish this section with the following remark.

\begin{Rem}[Generalization to modified policy iteration] 
Recall the generic value update rule 
\begin{equation}\label{eq:gen_update}
\bm{v}_{k+1} = \bm{v}_k + \bm{G}_k \big( \bo (\bm{v}_k) -  \bm{v}_k\big) , \quad k=0,1,\ldots,    
\end{equation}
with the gain matrix $\bm{G}_k = \bm{I}$ in the VI algorithm and $\bm{G}_k = (\bm{I} - \gamma \bm{P}_k)^{-1}$ in the PI algorithm. 
Also, note that $(\bm{I} - \gamma \bm{P}_k)^{-1} = \sum_{\ell=0}^{\infty} \gamma^\ell \bm{P}_k^\ell$ since $\rho(\bm{I} - \gamma \bm{P}_k) < 1$~\citep[Cor.~C.4]{puterman2014markov}. 
Inserting the truncated sum 
\[
\bm{G}_k = \sum_{\ell=0}^{L} \gamma^\ell \bm{P}_k^\ell, 
\]
with $L \in \{0,1,\ldots\}$ in the update rule~\eqref{eq:gen_update}, we derive the \emph{Modified PI} (\emph{MPI}) algorithm which converges linearly with rate $\gamma$ for any choice of $L$~\citep[Thm.~6.5.5]{puterman2014markov}. 
(Observe that $L=0$ and $L = \infty$ correspond to the standard VI and PI algorithms, respectively). 
The proposed rank-one modification can be in general combined with the MPI algorithm. 
Indeed, we have
\begin{align*}
 (\bm{I} - \gamma \bm{P}_k)^{-1} = \! \sum_{\ell=0}^{\infty} \gamma^\ell \bm{P}_k^\ell \!
   =  \sum_{\ell=0}^{L-1} \gamma^\ell \bm{P}_k^\ell + \gamma^L \bm{P}_k^L \sum_{\ell=0}^{\infty} \gamma^\ell \bm{P}_k^\ell  = \sum_{\ell=0}^{L-1} \gamma^\ell \bm{P}_k^\ell + \gamma^L \bm{P}_k^L (\bm{I} - \gamma \bm{P}_k)^{-1}.
\end{align*}
Then, by using the approximation $\wt{\bm{P}}_k = \e \bm{d}_k\tr$ in the matrix inversion on the right-hand side of the equation above, we derive the gain matrix of the \emph{rank-one MPI (R1-MPI)} algorithm to be 
\begin{align*}
  \bm{G}_k  =  \sum_{\ell=0}^{L-1} \gamma^\ell \bm{P}_k^\ell + \gamma^L \bm{P}_k^L (\bm{I} - \gamma \wt{\bm{P}}_k)^{-1}  = \sum_{\ell=0}^{L} \gamma^\ell \bm{P}_k^\ell + \frac{\gamma^{L+1}}{1-\gamma}\e \bm{d}_k\tr.
\end{align*}
Observe that R1-VI is now a special case of R1-MPI with $L=0$. 
\end{Rem}

\section{Rank-one Q-learning (R1-QL)}\label{sec:R1QL} 

In this section, we focus on the learning problem in which we have access to a generative model that provides us with samples of the MDP (as opposed to access to the true transition probability kernel of the MDP in the planning problem). 
To start, let us provide the PI update rule for the Q-function. 
The proof is similar to the proof of Lemma~\ref{lem:PI} and omitted. 
\begin{Lem}[Policy iteration for Q-function]\label{lem:PI-Q} 
Each iteration of the PI algorithm for the Q-function is given by 
\begin{equation}\label{eq:PI_Q}
    \bm{q}_{k+1} = \bm{q}_{k} + (\bm{I} - \gamma \ol{\bm{P}}_k)^{-1} \big(\ol{\bo}(\bm{q}_{k}) - \bm{q}_{k} \big),
\end{equation}
where~$\ol{\bm{P}}_k \Let \ol{\bm{P}}^{\pi^{\bm{q}_k}}$ is the state-action transition probability matrix of the MDP under the greed policy~$\pi^{\bm{q}_k}(s) \in \argmin_{a\in\set{A}} \bm{q}_k(s,a)$ for $s\in\set{S}$.
\end{Lem}

The idea is again to use the rank-one approximation~$\wt{\bm{P}}_k = \e \bm{d}_k\tr$ of the matrix~$\ol{\bm{P}}_k$ in the update rule~\eqref{eq:PI_Q}, where $\bm{d}_k$ is now the stationary distribution of $\ol{\bm{P}}_k$. 
This leads to the update rule 
\begin{equation*}
    \bm{q}_{k+1} = \ol{\bo}(\bm{q}_k) + \frac{\gamma}{1-\gamma} \big\langle \bm{d}_k, \ol{\bo}(\bm{q}_k) - \bm{q}_k \big\rangle \e,
\end{equation*}
at each iteration of the planning problem. 
Then, for the learning problem, considering the \emph{synchronous} update of all state-action pairs~$(s,a) \in \set{S}\times\set{A}$ at each iteration~$k$, we arrive at the \emph{rank-one Q-learning (R1-QL)} update rule 
\begin{equation}\label{eq:R1PL}
\begin{aligned}
    &\alpha_k = \frac{\gamma\lambda_k}{1-\gamma} \big\langle \wh{\bm{d}}_k ,  \sbo_k (\bm{q}_k) - \bm{q}_k \big\rangle, \\
    & \bm{q}_{k+1}  = (1-\lambda_k) \bm{q}_k + \lambda_k \sbo_k (\bm{q}_k) + \alpha_k \e,
\end{aligned}
\end{equation}
where $\lambda_k \geq 0$ are properly chosen step-sizes satisfying the Robbins–Monro conditions (e.g., $\lambda_k = 1/(k+1)$), 
and $\sbo_k$ is the empirical Bellman operator evaluated at iteration~$k$, i.e., $[\sbo_k (\bm{q}_k)](s,a) \Let [\sbo (\bm{q}_k,\wh{s}_k^{~+})](s,a)$ with $\wh{s}_k^{~+}\sim P(\cdot|s,a)$ for each $(s,a) \in\set{S}\times\set{A}$ -- the subscript $k$ in $\sbo_k$ denotes the dependence on the next-state sample $\wh{s}_k^{~+}$ generated at iteration~$k$. 

What remains to be addressed is computing the estimation~$\wh{\bm{d}}_k$ of the stationary distribution in~\eqref{eq:R1PL} using the samples. 
At each iteration $k$, define the sparse matrix~$\bm{F}_k \in \R^{|\set{S} \times \set{A}| \times |\set{S} \times \set{A}|} = \R^{(nm)\times(nm)}$ with exactly one nonzero entry equal to $1$ in each row~$(s,a) \in \set{S} \times \set{A}$ corresponding to the column $(s^+,a^+)$, where
\begin{equation*}
\begin{aligned}
&s^+ = \wh{s}_k^{~+} \sim P(\cdot|s,a),\\ 
&a^+ = \wh{a}_k^{~+} \in \argmin\limits_{a^+ \in \set{A}} \bm{q}_k(\wh{s}^{~+}_k,a^+).   
\end{aligned}
\end{equation*}
Observe that the matrix~$\bm{F}_k$ is a sampled version of the state-action transition probability matrix~$\ol{\bm{P}}_k$. 
Using this sample, we can form the stochastic approximation
\begin{equation*}
    \wh{\bm{P}}_k = (1-\lambda_k) \wh{\bm{P}}_{k-1} + \lambda_k \bm{F}_k.
\end{equation*}
for the state-action transition probability matrix. 
We note that the same approximation is used in the Zap Q-learning algorithm~\cite{devraj2017zap}. 
With this approximation in hand, we can again use the power method for finding the stationary distribution. 
In particular, with a single iteration of the power method initialized by the previous stationary distribution~$\wh{\bm{d}}_{k-1}$, we have
\begin{align*}
    \wh{\bm{d}}_{k} = \wh{\bm{P}}_k\tr  \wh{\bm{d}}_{k-1} 
    = (1-\lambda_k) \wh{\bm{P}}_{k-1}\tr \wh{\bm{d}}_{k-1} + \lambda_k \bm{F}_k\tr \wh{\bm{d}}_{k-1}.
\end{align*}
Now, using the approximation $\wh{\bm{d}}_{k-1} \approx \wh{\bm{P}}_{k-1}\tr \wh{\bm{d}}_{k-1}$ (i.e., assuming $\wh{\bm{d}}_{k-1}$ is a stationary distribution of $\wh{\bm{P}}_{k-1}$ which does not hold exactly since $\wh{\bm{d}}_{k-1}$ is only an approximation of a stationary distribution of $\wh{\bm{P}}_{k-1}$), we derive 
\begin{equation*}
\wh{\bm{d}}_{k} = (1-\lambda_k) \wh{\bm{d}}_{k-1} + \lambda_k \bm{F}_k\tr \wh{\bm{d}}_{k-1}.
\end{equation*}
We note that the vector $\bm{f} = \bm{F}_k\tr \wh{\bm{d}}_{k-1}$ can be computed using the following pseudo-code:
\begin{equation*}
    \begin{array}{l}
        \texttt{intialize } \bm{f} = \bm{0} \in \R^{|\set{S}\times\set{A}|} \\ 
        \texttt{for } (s,a)\in\set{S}\times\set{A}  \\ 
        \hspace{0.5cm} \wh{s}_k^{~+} \sim P(\cdot|s,a),\quad \wh{a}_k^{~+} \in \argmin\limits_{a^+ \in \set{A}} \bm{q}_k(\wh{s}^{~+}_k,a^+) \\
        \hspace{0.5cm} \bm{f}(\wh{s}_k^{~+},\wh{a}_k^{~+}) = \bm{f}(\wh{s}_k^{~+},\wh{a}_k^{~+}) + \wh{\bm{d}}_{k-1} (s,a)\\ 
        \texttt{endfor }
    \end{array}
\end{equation*}
Observe that the approximation $\wh{\bm{d}}_{k-1} \approx \wh{\bm{P}}_{k-1}\tr \wh{\bm{d}}_{k-1}$ significantly reduces the memory and time complexity of the algorithm since we do \emph{not} need to keep track of the estimates~$\wh{\bm{P}}_k$ of the state-action transition probability matrix and perform full matrix-vector multiplications for updating the estimates~$\wh{\bm{d}}_{k}$ of the stationary distribution. 

\begin{algorithm}[t]
\begin{small}
   \caption{Rank-One Q-Learning (R1-QL)} 
   \label{alg:R1-QL}
\begin{algorithmic}[1]
	\REQUIRE samples from transition kernel~$P: \set{S}\times\set{A}\ra \Delta(\set{S})$; cost function~$\bm{c} \in \R^{|\set{S}\times\set{A}|}$; discount factor~$\gamma \in (0,1)$;  
	\ENSURE optimal Q-function~$q\opt$
  	
	\vspace{.2cm}
  	
    \STATE initialize: $\bm{q}_0\in  \R^{|\set{S} \times \set{A}|},\; \wh{\bm{d}}_{-1} \in \Delta(\set{S}\times\set{A})$;

    \vspace{.1cm}
    
    \FOR{$k = 0, 1, 2, \ldots$} 
        \STATE $\lambda_k = 1/(k+1)$; 
        \STATE $\bm{f} = \z \in \R^{|\set{S} \times \set{A}|}$;
        \FOR{$(s,a) \in\set{S}\times\set{A}$}
            \STATE $\wh{s}_k^{~+}\sim P(\cdot|s,a)$;
            \vspace{.1cm}
            \STATE $\left\{\begin{array}{l}
                 \wh{a}_k^{~+} \in \argmin\limits_{a^+ \in \set{A}} \bm{q}_k(\wh{s}^{~+}_k,a^+), \\
                 \big[ \sbo_k (\bm{q}_k) \big] (s,a) = \bm{c}(s,a) + \gamma  \min\limits_{a^+ \in \set{A}} \bm{q}_k(\wh{s}^{~+}_k,a^+);
            \end{array}\right.$
            \vspace{.1cm}
            \STATE $\bm{f}(\wh{s}_k^{~+},\wh{a}_k^{~+}) = \bm{f}(\wh{s}_k^{~+},\wh{a}_k^{~+}) + \wh{\bm{d}}_{k-1} (s,a);$
        \ENDFOR
        \STATE $\wh{\bm{d}}_{k} = (1-\lambda_k) \wh{\bm{d}}_{k-1}  + \lambda_k\bm{f}$;  $\wh{\bm{d}}_k = \wh{\bm{d}}_{k} / \Vert\wh{\bm{d}}_{k}\Vert_1$; 
        \STATE $\alpha_k = \frac{\gamma\lambda_k}{1-\gamma} \big\langle \wh{\bm{d}}_k ,  \sbo_k (\bm{q}_k) - \bm{q}_k \big\rangle$;
        \STATE $\bm{q}_{k+1} = (1-\lambda_k) \bm{q}_k + \lambda_k \sbo_k (\bm{q}_k)  +  \alpha_k \e$; 
           
    \ENDFOR
   
\end{algorithmic}
\end{small}
\end{algorithm}

The complete description of the proposed R1-QL algorithm is provided in Algorithm~\ref{alg:R1-QL}. 
We again note that the normalization is only introduced to avoid the accumulation of numerical errors. 
The following result discusses the convergence of the proposed algorithm.

\begin{Thm}[Convergence of R1-QL]\label{thm:R1Ql-conv}
    The iterates~$\bm{q}_k$ of the R1-QL Algorithm~\ref{alg:R1-QL} converge to the optimal Q-function~$\bm{q}\opt = \ol{\bo}(\bm{q}\opt)$ almost surely with at least the same rate as QL.
\end{Thm}

Finally, we note that the per-iteration time complexity of the R1-QL Algorithm~\ref{alg:R1-QL} is the same as that of the synchronous QL algorithm, i.e., $\ord(nm^2)$.

\section{Numerical simulations}\label{sec:num}

In this section, we compare the performance of several planning and learning algorithms with our proposed methods. The experiments are conducted on Garnet \cite{archibald1995generation} and Graph MDPs \cite{devraj2017zap}, focusing on the Bellman errors~$\norm{\bo (\bm{v}_k) - \bm{v}_k}_\infty$ and~$\norm{\ol{\bo} (\bm{q}_k) - \bm{q}_k}_\infty$ and the value errors~$\norm{\bm{v}_k - \bm{v}\opt}_\infty$ and~$\norm{\bm{q}_k - \bm{q}\opt}_\infty$.
In all of our experiments, we run policy iteration (PI) until it converges to calculate the optimal values $\bm{v}\opt$ and $\bm{q}\opt$ for reference.
Garnet and Graph MDPs are particularly compelling for our empirical analysis as we can run PI to compute the optimal values, enabling us to measure value errors throughout the iterations.
The Garnet MDPs have a state size of $n = 200$, an action size of $m=5$, and randomly generated transition probabilities and costs with a branching factor of 10. 
For our numerical experiments, we consider 25 randomly generated instances of Garnet MDPs and report the three quantiles of the errors. 
For the Graph MDPs, we use the same configuration as described in~\cite{devraj2017zap}, providing a complementary benchmark to validate the effectiveness of our proposed method. 
In what follows, we report the result of our simulations for the planning and learning problems. 

\textbf{Planning Algorithms.} 
We compare several value iteration (VI) algorithms with the same per-iteration time complexity as our proposed R1-VI Algorithm~\ref{alg:R1-VI}. 
We also include PI for reference. 
We mainly focus on comparing R1-VI with accelerated VI methods, namely, Nesterov-VI~\cite{goyal2019first} and Anderson-VI~\cite{geist2018anderson}. 
In order to keep the time complexity the same, we use Anderson-VI with the memory parameter equal to one leading to a rank-one approximation of the Hessian matrix. 
The update rule of the accelerated VI algorithms is provided in Appendix~\ref{app:alg}.

\begin{figure*}[!htb]
    \centering
    \includegraphics[width=\textwidth]{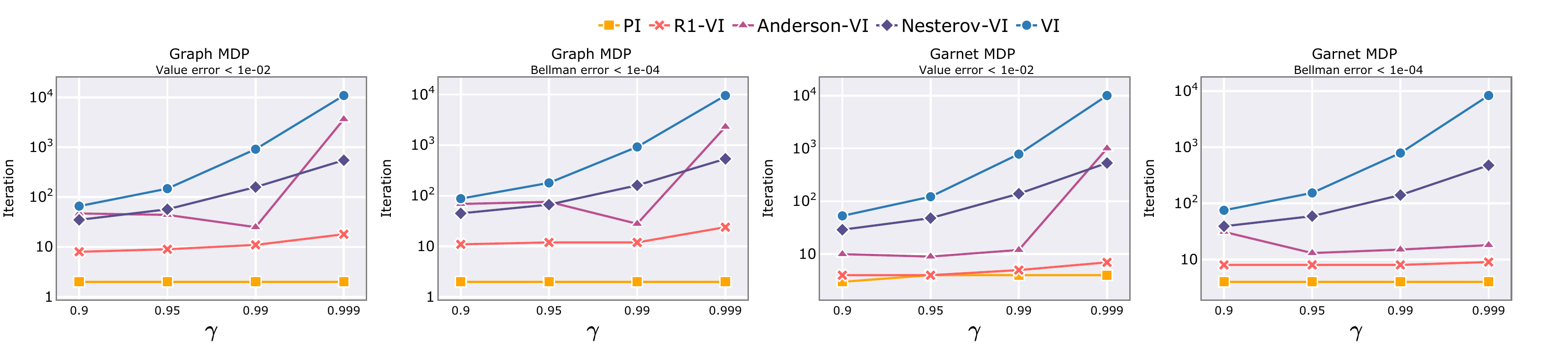}
    \caption{Planning algorithms -- the median number of iterations required for each algorithm to reach a fixed error threshold across four discount factors~$\gamma$ for the two MDPs.}
    \label{fig:planning-comparison-gamma}
\end{figure*}

\begin{figure*}[!htb]
    \centering
    \includegraphics[width=\textwidth]{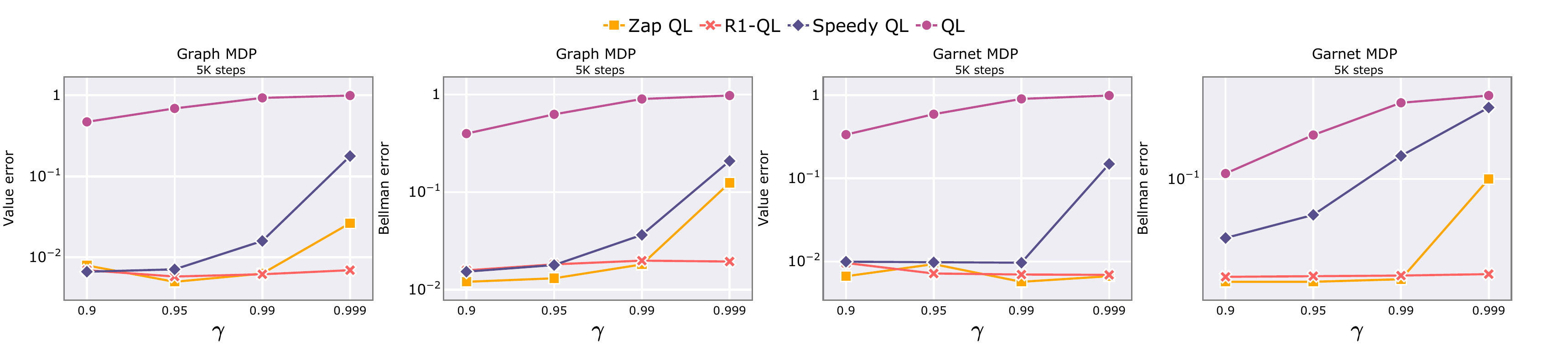}
    \caption{Learning algorithms -- the median error values achieved by each learning algorithm over the course of 5000 iterations across four discount factors~$\gamma$ for the two MDPs.}
    \label{fig:learning-comparison-gamma}
\end{figure*}

In Figure~\ref{fig:planning-comparison-gamma}, we report the median number of iterations required to reach a certain error threshold for each algorithm across both MDPs for four different values of the discount factors~$\gamma$. 
Our results indicate that the convergence performance of R1-VI is comparable with PI while maintaining the same per-iteration time complexity as VI. 
Additionally, R1-VI significantly outperforms the VI algorithm and its accelerated versions, particularly when the discount factor is close to 1. 
This observation can be partially explained by the fact that the proposed R1-VI algorithm forms an approximation of the inverse of the ``Hessian'', i.e.,~$(\bm{I} - \gamma \bm{P}_k)^{-1}$, by incorporating its largest eigenvalue~$\frac{1}{1-\gamma}$. 
A more detailed comparison of the algorithms is provided in Appendix~\ref{app:num_extra}.

\textbf{Learning Algorithms.}
In the learning experiments, we report the Bellman and value errors of the algorithms trained in a synchronous fashion, following the methodology outlined in~\cite{ghavamzadeh2011speedy}, ensuring a consistent evaluation.
In synchronous learning, in each iteration~$k$, a sample $\wh{s}^{~+} \sim P(\cdot \vert s, a)$ of the next state is generated for each state-action pair $(s, a)\in \set{S}\times\set{A}$ and the action-value function~$\bm{q}_k$ is updated in all state-action pairs; see the update rule in QL algorithm~\eqref{eq:QL} and R1-VI Algorithm~\ref{alg:R1-QL}.
All the learning algorithms use the same samples generated through the training.
Besides the proposed R1-QL Algorithm~\ref{alg:R1-QL}, we report the performance of Speedy QL~\cite{ghavamzadeh2011speedy}, Zap QL~\cite{devraj2017zap}, and the standard QL~\eqref{eq:QL} in Garnet and Graph MDPs for several discount factors; 
see Appendix~\ref{app:alg} for the update rule of Speedy QL and Zap QL. 
We run each algorithm using the same step-size schedule $(\lambda_k)_{k=0}^{\infty}$, namely, linearly decaying~$\lambda_k = 1/(1+k)$, to ensure a fair comparison.

Figure~\ref{fig:learning-comparison-gamma} shows the final error values after running each algorithm for 5000 iterations.
R1-QL achieves comparable or lower error values across both MDPs.
In contrast to the other algorithms, R1-QL consistently maintains a similar level of error values across various discount factors, particularly at higher discount factors -- a characteristic attributed to Policy Iteration (PI) algorithms.
It is worth mentioning that since both Zap QL and R1-QL estimate~$(\bm{I} - \gamma \ol{\bm{P}}_k)^{-1}$ using the samples, they behave more robustly against the increase in the discount factor~$\gamma$; see Figure~\ref{fig:learning-comparison-gamma}.
Regarding per iteration complexity, R1-QL has the same time and memory complexity as QL and Speedy QL. 
In contrast, Zap QL, due to the inherent full-rank matrix inversion, incurs higher time and memory complexity.
Albeit Zap QL can be efficiently implemented with lower complexity, it typically exhibits a higher computational cost in practice.
A comprehensive analysis of the error trajectories observed during training is provided in Appendix~\ref{app:num_extra}.

\section{Limitations and future research directions}\label{sec:lim}

We finish the paper by discussing some limitations of the proposed rank-one modification of the VI and QL algorithms along with some future research directions. 

Let us start by noting that the provided theoretical results in Theorems~\ref{thm:R1VI-conv} and~\ref{thm:R1Ql-conv} guarantee the convergence of the proposed algorithms with the same rate as the standard VI and QL algorithms. 
However, our numerical experiments with Garnet and Graph MDPs in Section~\ref{sec:num} show that the proposed algorithms have a faster convergence rate compared to standard VI and QL and their accelerated versions. 
This gap can be explained by the fact that our proof technique does not exploit that the vectors~$\bm{d}_k$ and~$\wh{\bm{d}}_k$ used in the update rule are specifically constructed to approximate the stationary distribution of the Markov chain induced by the greedy policy (see Appendices~\ref{app:proof_conv_R1VI} and~\ref{app:proof_conv_R1QL} for details). 
That is, the convergence of R1-VI and R1-QL algorithms is guaranteed for any choice of $\bm{d}_k \in \Delta(\set{S})$ and $\wh{\bm{d}}_k \in \Delta(\set{S}\times\set{A})$.
In fact, when the stationary distribution concentrates on a single state, as happens when there is an absorbing state with zero reward, the second term in the R1VI update rule~\eqref{eq:R1-VI_gen} vanishes. 
Appendix \ref{app:reducible} provides an empirical analysis in Gridworld \cite{sutton2018reinforcement}, which includes an absorbing state and thus violates the assumption of Lemma \ref{lem:unique_sol}.
Moreover, the provided proof of convergence shows that the greedy policies generated with respect to the iterates of R1-VI and R1-QL are the same as those for the standard VI and QL algorithms, respectively (see Lemmas~\ref{lem:simplified-R1-update} and \ref{lem:simplified-R1-updatemf}). 
In other words, the proposed algorithms do not affect the speed of convergence to the optimal policy compared to VI and QL. 
Nevertheless, at least in the case of R1-VI, the faster convergence in the value spaces leads to a faster termination of the algorithm for a given performance bound for the greedy policy.
In this regard, let us also note that the mismatch between convergence in value space and policy space also arises in other ``accelerated'' VI/QL algorithms; see Appendix~\ref{app:policy_performance}.

Second, the proposed algorithms heavily depend on the structure of the transition probability matrices~$ \bm{P}_k$ and~$\ol{\bm{P}}_k$ and their rank-one approximation using the corresponding stationary distributions. 
This dependence particularly hinders the application of the proposed algorithms to generic function approximation setups in solving the optimal control problem of MDPs with continuous state-action spaces. 
We note that a similar issue for the Zap Q-leaning algorithm~\cite{devraj2017zap} has been successfully addressed in~\cite{chen2020zap}. 

Third, we note that the proposed R1-QL algorithm~\ref{alg:R1-QL} is a \emph{synchronous} algorithm that updates \emph{all} state-action pairs~$(s,a) \in \set{S}\times\set{A}$ of the Q-function~$\bm{q}_k$ at each iteration~$k$. 
This algorithm can be modified in a standard fashion for the \emph{asynchronous} case. 
However, the provided convergence analysis can not be extended for the corresponding asynchronous algorithm in a straightforward manner.
Moreover, the straightforward asynchronous implementation of R1-QL leads to an $\ord(nm)$ per-iteration complexity for updating a \emph{single} component~$(s,a) \in \set{S}\times\set{A}$ at each iteration~$k$, which is higher than the~$\ord(m)$ per-iteration complexity of the standard asynchronous QL algorithm. 
We note that the Zap Q-leaning algorithm~\cite{devraj2017zap} also suffers from this issue. 
Addressing these issues requires a more involved analysis and modification of the proposed algorithm in the asynchronous case, which we leave for future research. 

Finally, the basic idea of the proposed algorithms can also be used for developing the rank-one modified version of the existing algorithms for the \emph{average} cost setting. 
For example, consider the PI algorithm that uses the \emph{relative VI} algorithm in the policy evaluation step for unichains \citep[Sec. 8.6.1]{puterman2014markov}. 
This algorithm can be characterized via the following update rule in the value space: For a fixed $s\in\mathcal{S}$,
\begin{align*}
   &\bm{v}_{k+1} = \bm{v}_k + \big((I-\bm{P}_k)(\bm{I} – e_{s} e_{s}^{\top})+\e e_{s}^{\top}\big)^{-1} \bm(\bo(\bm{v}_k) – \bm{v}_k\big),\\
   &\bm{v}_{k+1}( s) = 0,
\end{align*}
where $ e_{s}$ is the ${s}$-th unit vector and $\bo$ is now the \emph{undiscounted} Bellman operator. Now, observe that 
\begin{align*}
    \bm{G}_k = \big((\bm{I}-\bm{P}_k)(\bm{I} – e_{s} e_{s}^{\top})+\e e_{s}^{\top}\big)^{-1}
    = \big(\bm{I}-\bm{P}_k + (\bm{p}_k -  e_{s} +\e) e_{s}^{\top}\big)^{-1}
    = \big(\bm{I} - \e \bm{d}_k^{\top} + (\bm{p}_k -  e_{s}+\e) e_{s}^{\top}\big)^{-1},
\end{align*}
where $\bm{p}_k = \bm{P}_k(\cdot, s) $ is the $s$-th column of $\bm{P}_k $ and we used the approximation $\bm{P}_k \approx \e \bm{d}_k^{\top}$ in the last equality. 
The matrix inversion can then be handled efficiently using the Woodbury formula. 
However, the convergence of this algorithm and any possible improvement in the convergence rate when $\bm{d}_k$ is approximated via the power method requires further investigation.

\appendix

\section{Technical Proofs}\label{app:proofs}

\subsection{Proof of Lemma~\ref{lem:PI}}
 
We begin by providing two basic results on MDPs. 
First, recall that given a policy~$\pi$, the value~$\bm{v}^{\pi}$ of~$\pi$ solves the (\emph{Bellman consistency}) equation~\citep[Thm.~6.1.1]{puterman2014markov} 
\begin{equation*}
    \bm{v}^{\pi} (s) = \bm{c}^{\pi}(s) + \gamma  \ \EE_{s^+\sim \bm{P}^{\pi}(s,\cdot)} \left[ \bm{v}^{\pi}(s^+)\right],\quad  \forall s \in \set{S}.
\end{equation*}
Hence, we have
\begin{equation}\label{eq:eval}
    \bm{v}^{\pi} = (\bm{I} - \gamma \bm{P}^{\pi})^{-1} \bm{c}^{\pi}.
\end{equation}
Moreover, the definitions of the Bellman operator in~\eqref{eq:Bellman_operator} and of the greedy policy in~\eqref{eq:greedy_pol} imlpy that
\begin{equation}\label{eq:improve}
    \bo (\bm{v}) = \bm{c}^{\pi^{\bm{v}}} + \gamma \bm{P}^{\pi^{\bm{v}}}\bm{v}, \quad \forall \bm{v} \in \R^n.
\end{equation}
Recall~$\bm{P}_k \Let \bm{P}^{\pi^{\bm{v}_k}}$. Using these two results, we have at each iteration of the PI algorithm~\eqref{eq:PI},
\begin{align*}
    \bm{v}_{k+1} &= \bm{v}^{\pi_{k+1}} = \bm{v}^{\pi^{\bm{v}_k}} \overset{\eqref{eq:eval}}{=} (\bm{I} - \gamma \bm{P}_k)^{-1} \bm{c}^{\pi^{\bm{v}_k}} \overset{\eqref{eq:improve}}{=} (\bm{I} - \gamma \bm{P}_k)^{-1} \big( \bo (\bm{v}_k) -  \gamma \bm{P}_k\bm{v}_k \big)\\
     &= (\bm{I} - \gamma \bm{P}_k)^{-1} \big( \bo (\bm{v}_k) - \bm{v}_k + (\bm{I}- \gamma \bm{P}_k)\bm{v}_k \big)\\
     &= \bm{v}_{k} + (\bm{I} - \gamma \bm{P}_k)^{-1} \big(\bo(\bm{v}_{k}) - \bm{v}_{k} \big).
\end{align*}
This concludes the proof. 

\subsection{Proof of Lemma~\ref{lem:unique_sol}}

Let~$(\rho_i)_{i=1}^n$ be the eigenvalues of~$\bm{P}_k$ such that~$ |\rho_1| \geq |\rho_2| \geq \cdots \geq |\rho_n|$. 
Since~$\bm{P}_k$ is a row stochastic matrix, we have~$\rho_1 = 1$~\citep[Thm.~3.4.1]{gallager2011discrete}. 
Moreover, the assumption that~$\bm{P}_k$ is irreducible and aperiodic implies that~$|\rho_i| < 1$ for all~$i \neq 1$~\citep[Thm.~3.4.1]{gallager2011discrete}, that is,~$\rho_1 = 1$ is the unique eigenvalue of~$\bm{P}_k$ on the unit circle in the complex plane and all other eigenvalues lie inside the unit disc. 
The unique (up to scaling) right and left eigenvectors corresponding to~$\rho_1 = 1$ are the all-one vector~$\bm{1}$ and the stationary distribution~$\bm{d}_k$. 
From these results it follows that~$\wt{\bm{P}}_k = \e \bm{d}_k\tr$ is the unique solution of~\eqref{eq:QPI-R1 approx}. 

\subsection{Proof of Theorem~\ref{thm:R1VI-conv}}\label{app:proof_conv_R1VI}

For each iteration $k\geq 0$ of the R1-VI Algorithm~\ref{alg:R1-VI}, we have
\begin{align}
    \label{eq:R1-VI_1}
    \left\{
    \begin{aligned}
        &\bm{v}_{k+1} = \bo(\bm{v}_k) + \alpha_k \e \\
        &\alpha_k = \frac{\gamma}{1-\gamma} \big\langle \bm{d}_k, \bo(\bm{v}_k) - \bm{v}_k \big\rangle,
    \end{aligned}
    \right.
\end{align}
where $\bm{d}_k \in \Delta(\set{S})$ is an approximation of the stationary distribution of the MDP under the greedy policy with respect to~$\bm{v}_k$. 
We start with analyzing the effect of a constant shift in the argument of the Bellman operator.

\begin{Lem}\label{lem:linear_value}
    For all~$\alpha \in \R$ and~$\bm{v} \in \R^{n}$, we have
    \[
        \bo(\bm{v} + \alpha \e) = \bo(\bm{v}) + \gamma \alpha \e.
    \]
\end{Lem}
\begin{proof}
For each~$s \in \set{S}$, we have
\begin{align*}
    [\bo (\bm{v} + \alpha \e)](s) & = \min_{a \in\set{A}} \left\{ \bm{c}(s,a) + \gamma \sum_{s^+ \in \set{S}} P(s^+|s,a) \left[\bm{v} + \alpha \e \right] (s^+) \right\}\\
    & = \gamma \alpha + \min_{a \in\set{A}} \left\{ \bm{c}(s,a) + \gamma \sum_{s^+ \in \set{S}} P(s^+|s,a) \bm{v}  (s^+) \right\}\\
    & = \gamma \alpha + [\bo (\bm{v})](s).
\end{align*}
\end{proof}

We next use the preceding result to provide an alternative characterization of the iterates of R1-VI in~\eqref{eq:R1-VI_1}.

\begin{Lem}\label{lem:simplified-R1-update}
    For each~$k\geq 0$, the iterates in~\eqref{eq:R1-VI_1} are equivalently given by
    \[\left\{
    \begin{aligned}
        &\bm{v}_{k+1} = \bo^{(k+1)} (\bm{v}_0) + \beta_{k+1} \e \\
        &\beta_{k+1} = \gamma \beta_k + \alpha_k,\quad \text{\emph{with}}~ \beta_0 = 0.
    \end{aligned}
    \right.
    \]
\end{Lem}
\begin{proof}
(Proof by induction.) Consider~$k=0$ and observe that
\[
    \bm{v}_1 = \bo (\bm{v}_0) + \beta_1 \e,
\]
with $\beta_1 = \gamma \beta_0 + \alpha_0 = \alpha_0$. 
Next, for some~$k \geq 1$ and~$\beta_k \in \R$, assume $\bm{v}_k = \bo^{(k)} (\bm{v}_0) + \beta_k \e$. 
Then,
\begin{align*}
    \bm{v}_{k+1} & = \bo (\bm{v}_k) + \alpha_k \e
            = \bo \big(\bo^{(k)} (\bm{v}_0) + \beta_k \e \big) + \alpha_k \e 
            = \bo^{(k+1)} (\bm{v}_0) + \gamma \beta_k \e + \alpha_k \e,
\end{align*}
where the last equality follows from Lemma~\ref{lem:linear_value}. Therefore,
\begin{align*}
    \bm{v}_{k+1} & = \bo^{(k+1)} (\bm{v}_0) + (\gamma \beta_k + \alpha_k) \e 
            = \bo^{(k+1)} (\bm{v}_0) + \beta_{k+1} \e,
\end{align*}
which concludes the proof.   
\end{proof}

We now employ the preceding lemmas to provide an alternative characterization of the constant shifts~$\alpha_k$ in the R1-VI updates in~\eqref{eq:R1-VI_1}. 

\begin{Lem}\label{lem:stepsize-update}
    For each $k\geq 0$, one has
    \[
        \alpha_k = \frac{\gamma}{1 - \gamma} \left\langle \bm{d}_k, \bo^{(k+1)} (\bm{v}_0) - \bo^{(k)} (\bm{v}_0)  \right\rangle - \gamma \beta_k.
    \]
\end{Lem}
\begin{proof}
From Lemma~\ref{lem:simplified-R1-update}, we have
\[
    \bm{v}_k = \bo^{(k)}(\bm{v}_0) + \beta_k \e, \quad  k=0,1,\ldots .
\]
(Notice that for~$k=0$, the preceding equation simply implies~$\bm{v}_0=\bm{v}_0$ since~$\bo^{(0)}$ is the identity operator and~$\beta_0=0$.) 
Then, from Lemma~\ref{lem:linear_value}, it follows that
\[
    \bo (\bm{v}_k) = \bo \big(\bo^{(k)} (\bm{v}_0) + \beta_k  \e \big) = \bo^{(k+1)} (\bm{v}_0) + \gamma \beta_k  \e.
\]
Thus,
\[
    \bo(\bm{v}_k) - \bm{v}_k = \bo^{(k+1)} (\bm{v}_0) - \bo^{(k)} (\bm{v}_0)  - (1 - \gamma) \beta_k  \e
\]
and
\begin{align*}
    \big\langle \bm{d}_k , \bo(\bm{v}_k) - \bm{v}_k  \big\rangle 
    & = \big\langle \bm{d}_k,\bo^{(k+1)} (\bm{v}_0) - \bo^{(k)} (\bm{v}_0)  \big\rangle - (1 - \gamma) \beta_k \langle \bm{d}_k,\e \rangle \\
    & = \big\langle \bm{d}_k,\bo^{(k+1)} (\bm{v}_0) - \bo^{(k)} (\bm{v}_0) \big\rangle - (1 - \gamma) \beta_k, 
\end{align*}
where we used~$\bm{d}_k \in \Delta (\set{S})$ (i.e.,~$\langle \bm{d}_k ,  \e \rangle=1$) in the second equality above. 
Therefore, for $\alpha_k$ in~\eqref{eq:R1-VI_1}, we have
\[
    \alpha_k = \frac{\gamma}{1 - \gamma} \left\langle \bm{d}_k, \bo^{(k+1)} (\bm{v}_0) - \bo^{(k)} (\bm{v}_0) \right\rangle - \gamma \beta_k.
\]
This completes the proof.     
\end{proof}

Now, observe that plugging in $\alpha_k$ from Lemma~\ref{lem:stepsize-update} in the update rule of Lemma~\ref{lem:simplified-R1-update} leads to 
\begin{equation}\label{eq:R1-VI_spec}
        \bm{v}_{k+1} = \bo^{(k+1)}(\bm{v}_0) + \frac{\gamma}{1 - \gamma} \big\langle \bm{d}_k ,  \bo^{(k+1)}(\bm{v}_0) -\bo^{(k)}(\bm{v}_0) \big\rangle  \e, \quad k=0,1,\ldots.
\end{equation}
Then, using the fact that $\bm{v}\opt = \bo (\bm{v}\opt)$ and the Bellman operator is a $\gamma$-contraction in the $\infty$-norm, we have
\begin{align*}
    \norm{\bm{v}_{k+1} - \bm{v}\opt}_{\infty} & = \norm{  \bo^{(k+1)}(\bm{v}_0) - \bo (\bm{v}\opt) + \frac{\gamma}{1 - \gamma} \big\langle \bm{d}_k ,  \bo^{(k+1)}(\bm{v}_0) -\bo^{(k)}(\bm{v}_0) \big\rangle  \e }_{\infty} \\
    & \leq   \norm{\bo^{(k+1)}(\bm{v}_0) - \bo (\bm{v}\opt)}_{\infty} +  \frac{\gamma}{1 - \gamma} \left| \big\langle \bm{d}_k ,  \bo^{(k+1)}(\bm{v}_0) -\bo^{(k)}(\bm{v}_0) \big\rangle \right| \\
    & \leq  \gamma^{k+1} \norm{\bm{v}_0 - \bm{v}\opt}_{\infty} +  \frac{\gamma}{1 - \gamma} \norm{\bo^{(k+1)}(\bm{v}_0) -\bo^{(k)}(\bm{v}_0) }_{\infty} \\
    & \leq  \gamma^{k+1} \norm{\bm{v}_0 - \bm{v}\opt}_{\infty} +  \frac{\gamma^{k+1}}{1 - \gamma} \norm{\bo(\bm{v}_0) -\bm{v}_0}_{\infty}\\
    & \leq  \gamma^{k+1} \big( \norm{\bm{v}_0 - \bm{v}\opt}_{\infty} +  \frac{1}{1 - \gamma} \norm{\bo(\bm{v}_0) -\bm{v}_0}_{\infty} \big).
\end{align*}
That is, $\bm{v}_k \to \bm{v}\opt$ as~$k\to\infty$ linearly with rate $\gamma$.

\subsection{Proof of Theorem~\ref{thm:R1Ql-conv}}\label{app:proof_conv_R1QL}

Let us begin with recalling the definition of the empirical Bellman operator
\begin{align}\label{eq:sampled_bellman}
    \big[ \sbo_k (\bm{q}) \big] (s,a) := \bm{c}(s,a) + \gamma \min_{a^+ \in \set{A}} \bm{q}(\wh{s}^{~+}_k,a^+),
\end{align}
where~$\wh{s}_k^{~+} \sim P(\cdot | s,a)$, that is to say~$\wh{s}^{~+}$ is sampled according to the law~$P(\cdot | s,a)$ at iteration~$k$. 
From this definition, it immediately follows that 
\begin{equation}\label{lem:sam_Bell_linearity}
\sbo_k (\bm{q} + \alpha \e) = \sbo_k (\bm{q}) + \gamma \alpha \e, \quad \forall\ \alpha \in \R,\ \bm{q} \in \R^{|\set{S}\times\set{A}|}.   
\end{equation}
Also, recall that each iteration $k\geq 0$ of the R1-QL Algorithm~\ref{alg:R1-QL} reads as
\begin{align}\label{eq:R1-QL_mfnew}
    \left\{
    \begin{aligned}
        &\bm{q}_{k+1} = (1 - \lambda_k) \bm{q}_k + \lambda_k \sbo_k (\bm{q}_k) + \alpha_k \e \\
        &\alpha_k = \lambda_k \Big( \frac{\gamma}{1-\gamma} \Big) \big\langle \wh{\bm{d}}_k , \sbo_k (\bm{q}_k) -  \bm{q}_k \big\rangle,
    \end{aligned}
    \right.
\end{align}
where $\wh{\bm{d}}_k \in \Delta(\set{S}\times\set{A})$ is an estimation of the stationary distribution of the state-action transition probability matrix of the MDP under the greedy policy with respect to~$\bm{q}_k$. 
Let us also consider the standard QL iterates 
\[
    \bm{q}^{\text{QL}}_{k+1} = (1-\lambda_k) \bm{q}^{\text{QL}}_{k} + \lambda_k \sbo_k (\bm{q}^{\text{QL}}_{k}), \quad k=0,1,\ldots,
\]
with the same initialization~$\bm{q}^{\text{QL}}_0 = \bm{q}_0$ and empirical Bellman operator~$\sbo_k(\cdot)$ for all $k$ as the R1-QL algorithm~\eqref{eq:R1-QL_mfnew}. 

The first result concerns an alternative characterization of the iterates in~\eqref{eq:R1-QL_mfnew}.

\begin{Lem}\label{lem:simplified-R1-updatemf}
    For each $k\geq 0$, the iterates of R1-QL algorithm~\eqref{eq:R1-QL_mfnew} equivalently read as
    \begin{equation}\label{eq:simplified-R1-updatemf}
       \left\{
    \begin{aligned}
        &\bm{q}_{k+1} = \bm{q}^{\emph{\text{QL}}}_{k+1} + \beta_{k+1} \e \\
        &\beta_{k+1} = (1 - \lambda_k) \beta_k + \gamma \lambda_k \beta_k + \alpha_k,\quad \text{\emph{with}}~ \beta_0 = 0.
    \end{aligned}
    \right. 
    \end{equation}
\end{Lem}
\begin{proof}
 (Proof by induction) For~$k=0$, since $\bm{q}^{\text{QL}}_0 = \bm{q}_0$, we can write
\begin{align*}
    \bm{q}_1 & = (1-\lambda_0) \bm{q}_0 + \lambda_0 \sbo_k (\bm{q}_0)+ \alpha_0 \e \\
        & = (1-\lambda_0) \bm{q}^{\text{QL}}_0 + \lambda_0 \sbo_k (\bm{q}^{\text{QL}}_0)  + \alpha_0 \e \\
        & = \bm{q}^{\text{QL}}_1 + \beta_1 \e,
\end{align*}
where~$\beta_1 = \alpha_0$. 
Assume next~$\bm{q}_k = \bm{q}^{\text{QL}}_k + \beta_k \e$ for some~$k \geq 0$. Then, it follows that
\begin{align*}
    \bm{q}_{k+1} & = (1 - \lambda_k) \bm{q}_k + \lambda_k \sbo_k (\bm{q}_k) + \alpha_k \e \\
            & = (1 - \lambda_k) ( \bm{q}^{\text{QL}}_k + \beta_k \e ) + \lambda_k \sbo_k ( \bm{q}^{\text{QL}}_k + \beta_k \e ) + \alpha_k \e \\
            & \overset{\eqref{lem:sam_Bell_linearity}}{=} (1 - \lambda_k) ( \bm{q}^{\text{QL}}_k + \beta_k \e ) + \lambda_k \big( \sbo_k ( \bm{q}^{\text{QL}}_k ) + \gamma \beta_k \e \big) + \alpha_k \e \\
            & = (1 - \lambda_k) \bm{q}^{\text{QL}}_k + \lambda_k \sbo_k ( \bm{q}^{\text{QL}}_k ) + \big( (1 - \lambda_k) \beta_k + \gamma \lambda_k \beta_k + \alpha_k \big) \e \\
            & = \bm{q}^{\text{QL}}_{k+1} +\beta_{k+1} \e.
\end{align*}
This concludes the proof.   
\end{proof}
We next provide a useful characterization of the constant shifts~$\alpha_k$ in the R1-QL update rule~\eqref{eq:R1-QL_mfnew}.

\begin{Lem}\label{lem:stepsize-updatemf}
    For each $k\geq 0$, one has
    \[
        \alpha_k = \lambda_k \Big( \frac{\gamma}{1 - \gamma} \Big) \big\langle \wh{\bm{d}}_k ,  \sbo_k ( \bm{q}^{\text{\emph{QL}}}_k) - \bm{q}^{\text{\emph{QL}}}_k  \big\rangle - \gamma \lambda_k \beta_k.
    \]
\end{Lem}
\begin{proof}
From Lemma~\ref{lem:simplified-R1-updatemf} and since $\bm{q}_0 = \bm{q}^{\text{QL}}_0$ and $\beta_0 = 0$, we have
\begin{align*}
    \bm{q}_k = \bm{q}^{\text{QL}}_k + \beta_k \e,\quad  k=0,1, \ldots.
\end{align*}
Hence, we can use \eqref{lem:sam_Bell_linearity} to write
\begin{align*}
    \sbo_k (\bm{q}_k) & = \sbo_k ( \bm{q}^{\text{QL}}_k + \beta_k \e ) 
                         = \sbo_k ( \bm{q}^{\text{QL}}_k ) + \gamma \beta_k \e. 
\end{align*}
As a result,
\begin{align*}
   \sbo_k (\bm{q}_k) - \bm{q}_k  & =  \sbo_k ( \bm{q}^{\text{QL}}_k ) + \gamma \beta_k \e - \bm{q}^{\text{QL}}_k - \beta_k \e  \\
                              & = \sbo_k ( \bm{q}^{\text{QL}}_k ) - \bm{q}^{\text{QL}}_k - (1 - \gamma) \beta_k \e,
\end{align*}
and
\begin{align*}
    \big\langle \wh{\bm{d}}_k ,\sbo_k (\bm{q}_k) - \bm{q}_k \big\rangle 
        & = \big\langle \wh{\bm{d}}_k , \sbo_k ( \bm{q}^{\text{QL}}_k ) - \bm{q}^{\text{QL}}_k  \big\rangle - (1 - \gamma) \beta_k \langle \wh{\bm{d}}_k , \e \rangle \\
        & = \big\langle \wh{\bm{d}}_k , \sbo_k ( \bm{q}^{\text{QL}}_k ) - \bm{q}^{\text{QL}}_k \big\rangle - (1 - \gamma) \beta_k,
\end{align*}
where we use the identity~$\langle \wh{\bm{d}}_k , \e \rangle=1$ in the second line since~$\wh{\bm{d}}_k \in \Delta(\set{S} \times \set{A})$. 
Recalling the definition of~$\alpha_k$ in~\eqref{eq:R1-QL_mfnew}, one thus have 
\begin{align*}
    \alpha_k = \lambda_k \Big( \frac{\gamma}{1 - \gamma} \Big) \big\langle \wh{\bm{d}}_k , \sbo_k ( \bm{q}^{\text{QL}}_k ) - \bm{q}^{\text{QL}}_k  \big\rangle - \gamma \lambda_k \beta_k.
\end{align*}    
\end{proof}
Plugging the expression for $\alpha_k$ from Lemma~\ref{lem:stepsize-updatemf} into the update rule of Lemma~\ref{lem:simplified-R1-updatemf}, we derive the iteration 
\begin{equation}\label{eq:beta}
    \beta_{k+1} = (1 - \lambda_k) \beta_k + \lambda_k \Big( \frac{\gamma}{1-\gamma} \Big) \big\langle \wh{\bm{d}}_k , \sbo_k ( \bm{q}^{\text{QL}}_k ) -\bm{q}^{\text{QL}}_k \big\rangle, \quad k = 0,1,\ldots, 
\end{equation}
initialized by $\beta_0 = 0$. 
Next, using the convergence of QL, we can show the convergence of $\beta_k$:  
\begin{Lem}
    The iterates $\beta_k$ in~\eqref{eq:beta} converge to zero almost surely. 
\end{Lem}
\begin{proof}

Define $\beta_{1,0} = \beta_{2,0} = 0$ and consider the iterations
\begin{align}
    \beta_{1,k+1} &= (1 - \lambda_k) \beta_{1,k} + \lambda_k \Big( \frac{\gamma}{1-\gamma} \Big) \big\langle \wh{\bm{d}}_k , \sbo_k ( \bm{q}^{\text{QL}}_k ) -\ol{\bo}(\bm{q}^{\text{QL}}_k) \big\rangle, \label{eq:beta_1} \\
    \beta_{2,k+1} &= (1 - \lambda_k) \beta_{2,k} + \lambda_k \Big( \frac{\gamma}{1-\gamma} \Big) \big\langle \wh{\bm{d}}_k , \ol{\bo}_k ( \bm{q}^{\text{QL}}_k ) -\bm{q}^{\text{QL}}_k \big\rangle, \label{eq:beta_2}
\end{align}
for $k = 0,1,\ldots$, so that $\beta_{k} = \beta_{1,k} + \beta_{2,k}$ for all $k\geq 0$. 
In what follows, we use the fact that the QL iterates $\bm{q}^{\text{QL}}_k$ are bounded and converge to $\bm{q}\opt = \ol{\bo}(\bm{q}\opt)$ almost surely~\citep[Thm.~4]{tsitsiklis1994asynchronous}. 
First, observe that the iteration~\eqref{eq:beta_1} converges to zero almost surely using~\citep[Lem.~1]{tsitsiklis1994asynchronous} and the fact that (see also ~\citep[Sec.~7]{tsitsiklis1994asynchronous}) 
\begin{align*}
    &\EE_{\wh{s}_k^+} \left[\big\langle \wh{\bm{d}}_k , \sbo_k ( \bm{q}^{\text{QL}}_k ) -\ol{\bo}(\bm{q}^{\text{QL}}_k) \big\rangle\right] = \big\langle \wh{\bm{d}}_k , \EE_{\wh{s}_k^+} \left[\sbo_k ( \bm{q}^{\text{QL}}_k ) -\ol{\bo}(\bm{q}^{\text{QL}}_k)\right] \big\rangle = 0,\\ 
    &\EE_{\wh{s}_k^+} \left[\big\langle \wh{\bm{d}}_k , \sbo_k ( \bm{q}^{\text{QL}}_k ) -\ol{\bo}(\bm{q}^{\text{QL}}_k) \big\rangle^2\right] \leq \EE_{\wh{s}_k^+} \left[\norm{\sbo_k ( \bm{q}^{\text{QL}}_k ) -\ol{\bo}(\bm{q}^{\text{QL}}_k)}_{\infty}^2\right] \leq \norm{\bm{q}^{\text{QL}}_k}_{\infty}^2.
\end{align*}
The iteration~\eqref{eq:beta_2} also converges to zero almost surely since 
\[
\beta_{2,k+1} = \Big(\frac{\gamma}{1-\gamma}\Big)  \sum_{(s,a)\in\set{S}\times\set{A}} \wh{\bm{d}}_k(s,a) \left\{ \frac{1}{k+1} \sum_{\ell=0}^{k} \big( [\ol{\bo}_{\ell} ( \bm{q}^{\text{QL}}_{\ell} )-\bm{q}^{\text{QL}}_{\ell}](s,a) \big)\right\},
\]
is a scaled, weighted average of Cesaro means of the sequences~$[\ol{\bo}_k ( \bm{q}^{\text{QL}}_k )-\bm{q}^{\text{QL}}_k](s,a)$ that converge to zero almost surely for each $(s,a)\in\set{S}\times\set{A}$ (recall that $\lambda_k = 1/(k+1)$ and $\bm{q}^{\text{QL}}_k\ra \bm{q}\opt = \ol{\bo}(\bm{q}\opt)$ almost surely).
\end{proof}

Finally, recall the characterization~$\bm{q}_{k} = \bm{q}^{\text{QL}}_{k} + \beta_{k} \e$ in Lemma~\ref{lem:simplified-R1-updatemf} and observe that $\bm{q}^{\text{QL}}_k\ra \bm{q}\opt$ and $\beta_{k} \ra 0$ almost surely. 
Therefore, $\bm{q}_{k}\ra\bm{q}\opt$ almost surely. 

\section{On numerical experiments}\label{app:num}

\subsection{Algorithms}\label{app:alg}

Below, we provide the update rules of the algorithms we employed in our numerical experiments. We adapt the update rules provided by~\cite{kolarijani2023optimization} in our implementations for the planning algorithms. 

\begin{itemize}
\item \textbf{Nesterov VI algorithm~\cite{goyal2019first}}: 
\begin{align*}
    \begin{aligned}
        &\bm{z}_{k} = \bm{v}_k + \frac{1 - \sqrt{1 - \gamma^2}}{\gamma} (\bm{v}_k - \bm{v}_{k-1}), \\
        &\bm{v}_{k+1} = \bm{z}_{k} + \frac{1}{1 + \gamma} (\bo(\bm{z}_{k}) - \bm{z}_{k}).
    \end{aligned}
\end{align*}

\item \textbf{Anderson VI algorithm \cite{geist2018anderson}}: (The following update rule is for Anderson acceleration with memory equal to $1$ which corresponds to a rank-one approximation of the Hessian.) 

\begin{align*}
    \begin{aligned}
        &\bm{z}_{k} = \bm{v}_k - \bm{v}_{k-1}, \\
        &\bm{z'}_{k} = \bo(\bm{v}_k) - \bo(\bm{v}_{k-1}), \\
        &\delta_k = \begin{cases}
            0, & \bm{z}_k^\top (\bm{z}_k - \bm{z}_k') = 0, \\
            \frac{\bm{z}_k^\top \big(\bm{v}_k - \bo(\bm{v}_k)\big)}{\bm{z}_k^\top (\bm{z}_k - \bm{z}_k')}, & \text{otherwise},
        \end{cases} \\
        &\bm{v}_{k+1} = (1 - \delta_k)\bo(\bm{v}_{k}) + \delta_k \bo(\bm{v}_{k - 1}).
    \end{aligned}
\end{align*}

\item \textbf{Speedy QL algorithm \cite{ghavamzadeh2011speedy}}: (The following update rule is the synchronous implementation of Speedy QL.)

\begin{align*}
    \begin{aligned}
        & \texttt{for } (s, a) \in \set{S} \times \set{A} \\
        & \qquad \wh{s}^{~+} \sim P(\cdot \vert s, a), \\
        & \qquad \bm{z}_k(s, a) = \bm{c}(s, a) + \gamma \min_{a^+ \in \set{A}} \bm{q}_k(\wh{s}^{~+} ,a^+), \\
        & \qquad \bm{z}_k'(s, a) = \bm{c}(s, a) + \gamma \min_{a^+ \in \set{A}} \bm{q}_{k-1}(\wh{s}^{~+}, a^+), \\
        & \texttt{endfor } \\ 
        & \qquad \bm{q}_{k+1} = \bm{q}_{k} + \frac{1}{1 + k} (\bm{z}_k' - \bm{q}_{k}) + \frac{k}{1 + k} (\bm{z}_k - \bm{z}_k').
    \end{aligned}
\end{align*}

\item \textbf{Zap QL algorithm \cite{devraj2017zap}}: (The following update rule is also the synchronous implementation of Zap QL without eligibility trace. 
The matrix~$\bm{F}_k \in \R^{|\set{S} \times \set{A}| \times |\set{S} \times \set{A}|}$ below denotes the sampled transition matrix at iteration~$k$.)

\begin{align*}
    \begin{aligned}
        & \bm{F}_k = \z \in \R^{|\set{S} \times \set{A}| \times |\set{S} \times \set{A}|}, \\
        & \texttt{for } (s, a) \in \set{S} \times \set{A} \\
        & \qquad \wh{s}^{~+} \sim P(\cdot \vert s, a), \\
        & \qquad \wh{a}^{~+} = \argmin_{a^+ \in \set{A}} \bm{q}_k(\wh{s}^{~+} ,a^+), \\
        & \qquad \big[ \sbo_k (\bm{q}_k) \big] (s,a) = \bm{c}(s, a) + \gamma \min_{a^+ \in \set{A}} \bm{q}_k(\wh{s}^{~+} ,a^+) - \bm{q}_{k}(s, a), \\
        & \qquad \bm{F}_k\big((s,a),(\wh{s}^{~+},\wh{a}^{~+})\big) = 1, \\
        & \texttt{endfor } \\ 
        & \wh{\bm{P}}_k = \wh{\bm{P}}_{k-1} + \frac{1}{2 + k} (\bm{F}_k - \wh{\bm{P}}_{k-1}),\\
        &\bm{q}_{k+1} = \bm{q}_{k} + \frac{1}{1 + k} (\bm{I} - \gamma \wh{\bm{P}}_k)^{-1} \big( \sbo_k (\bm{q}_k) - \bm{q}_{k} \big).
    \end{aligned}
\end{align*}

\end{itemize}

\subsection{Extended numerical analysis}\label{app:num_extra}

In this appendix, we provide the Bellman and value errors observed throughout the iterations of planning and learning algorithms. 
We run each planning algorithm until the error thresholds in Table \ref{tab:appx-error-treshold} are achieved.

\begin{table}[ht]
\centering
\caption{Error thresholds for planning algorithms.}
\label{tab:appx-error-treshold}
\begin{tabular}{@{}llccccc@{}}
\toprule
MDP & Error & $\gamma=0.9$ & $\gamma=0.9$ & $\gamma=0.9$ & $\gamma=0.9$ \\ \midrule

Garnet & value  & $10^{-5}$ & $10^{-4}$ & $10^{-4}$ & $10^{-2}$ \\
& Bellman  & $10^{-5}$ & $10^{-5}$ & $10^{-5}$ & $10^{-4}$\\
Graph & value  & $10^{-5}$ & $10^{-4}$ & $10^{-3}$ & $10^{-2}$ \\
& Bellman  & $10^{-5}$ & $10^{-5}$ & $10^{-5}$ & $10^{-4}$\\

\bottomrule
\end{tabular}
\end{table}

\begin{figure*}[!htb]
    \centering
    \includegraphics[width=\textwidth]{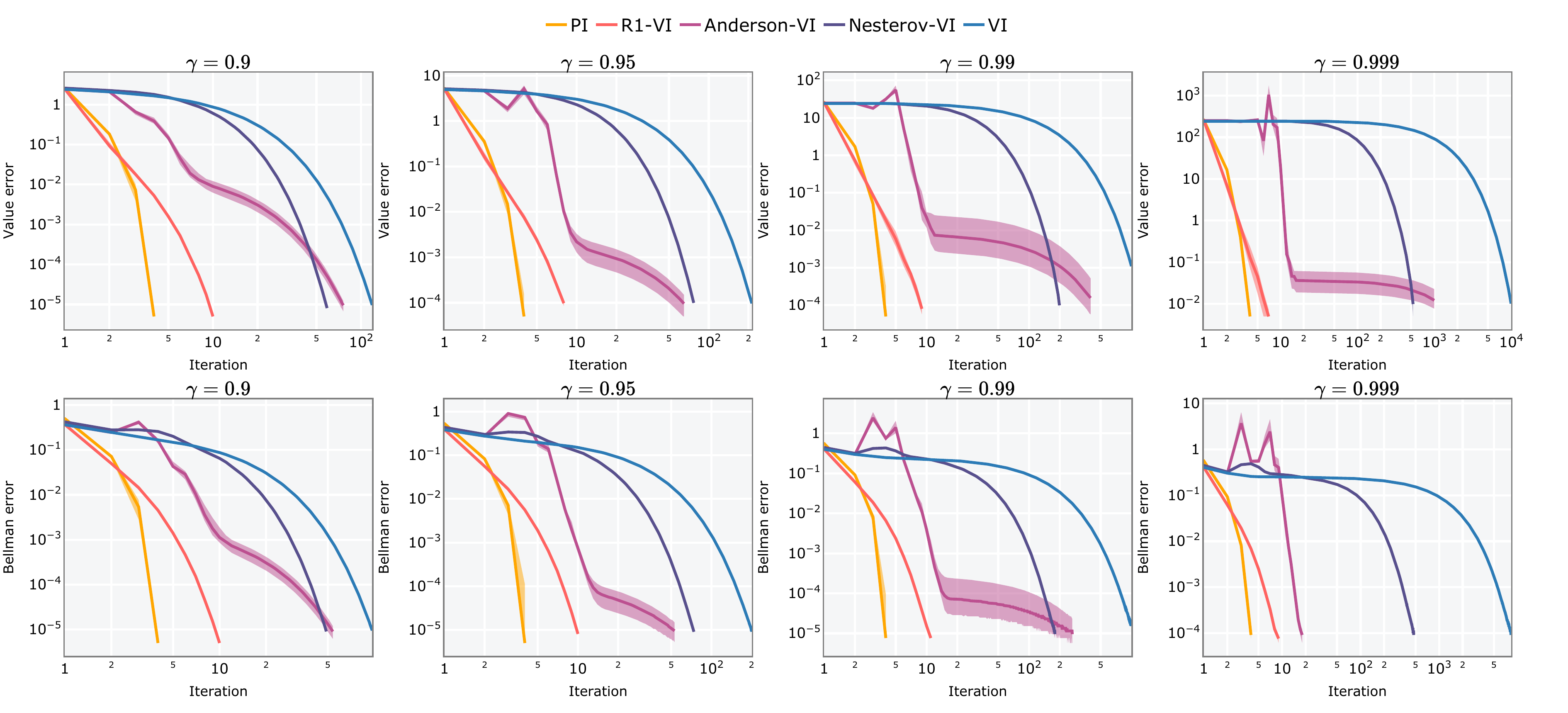}%
    \caption{Comparison of the planning algorithms in Garnet MDP with various $\gamma$ values.}
    \label{fig:appx-garnet-plan-curves}
\end{figure*}

We consider four different values for discount factor $\gamma$ across 25 realizations of the Garnet MDP. 
Note that the planning algorithms are deterministic in nature, hence, the only source of variation in the errors is due to the random realization of the Garnet MDPs. 
Figure \ref{fig:appx-garnet-plan-curves} shows the range of error values observed within the span of the iterations of the planning algorithms. 
The solid curve represents the median error values, while the shaded region around the curve indicates the errors between the first and the third quantiles. 
We follow the same style of presentation in the other figures in this section. 

\begin{figure*}[!htb]
    \centering
    \includegraphics[width=\textwidth]{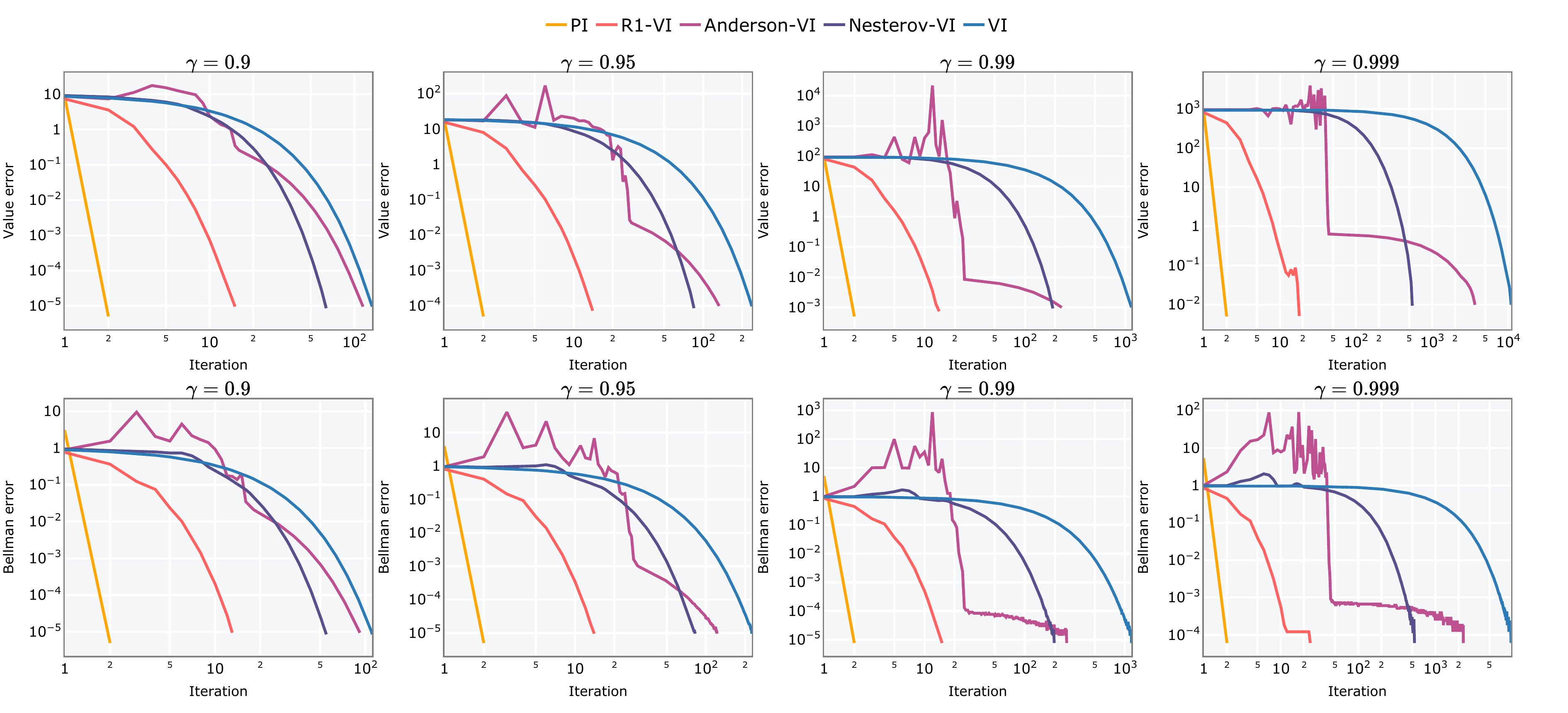}%
    \caption{Comparison of the planning algorithms in Graph MDP with various $\gamma$ values.}
    \label{fig:appx-graph-plan-curves}
\end{figure*}

\begin{figure*}[!htb]
    \centering
    \includegraphics[width=\textwidth]{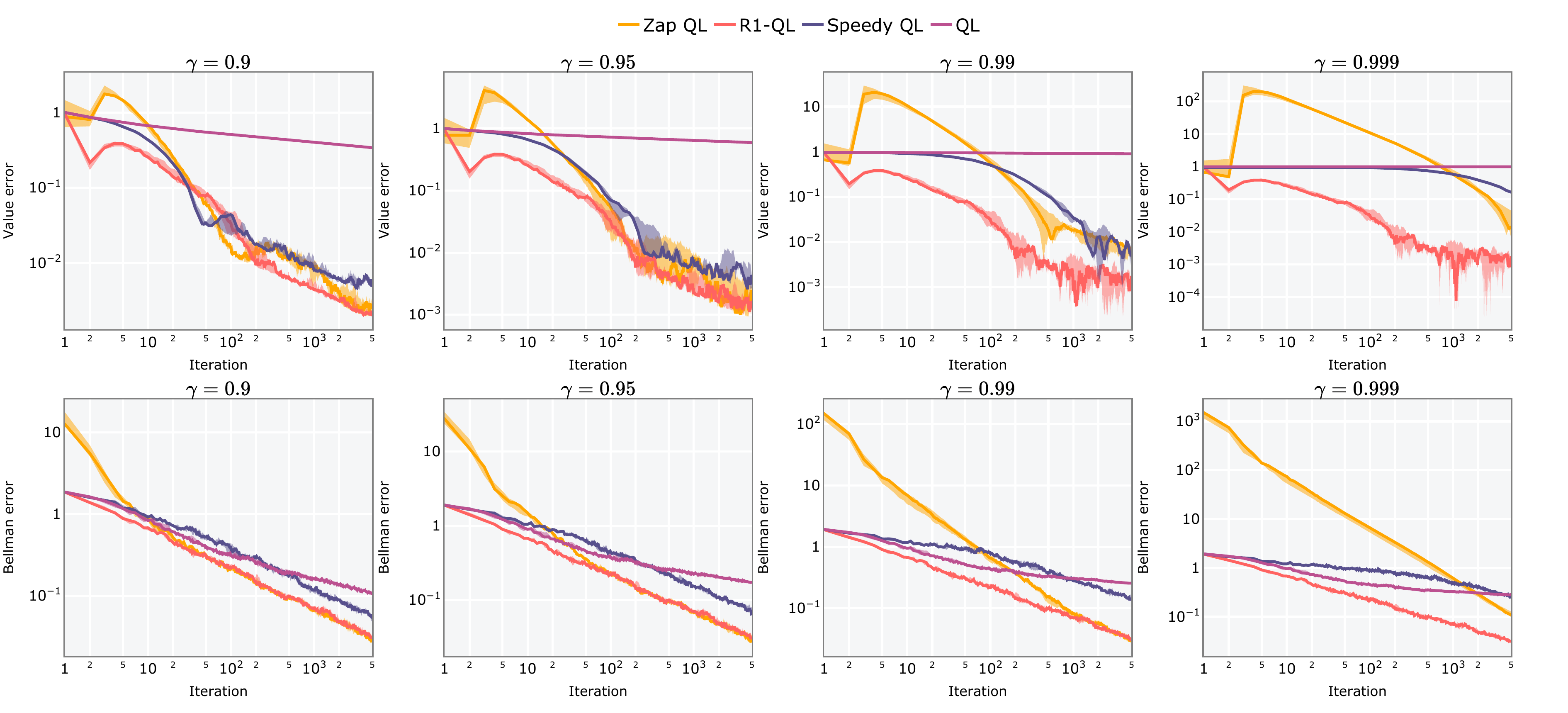}%
    \caption{Comparison of the learning algorithms in Garnet MDP with various $\gamma$ values.}
    \label{fig:appx-garnet-learn-curves}
\end{figure*}

\begin{figure*}[!htb]
    \centering
    \includegraphics[width=\textwidth]{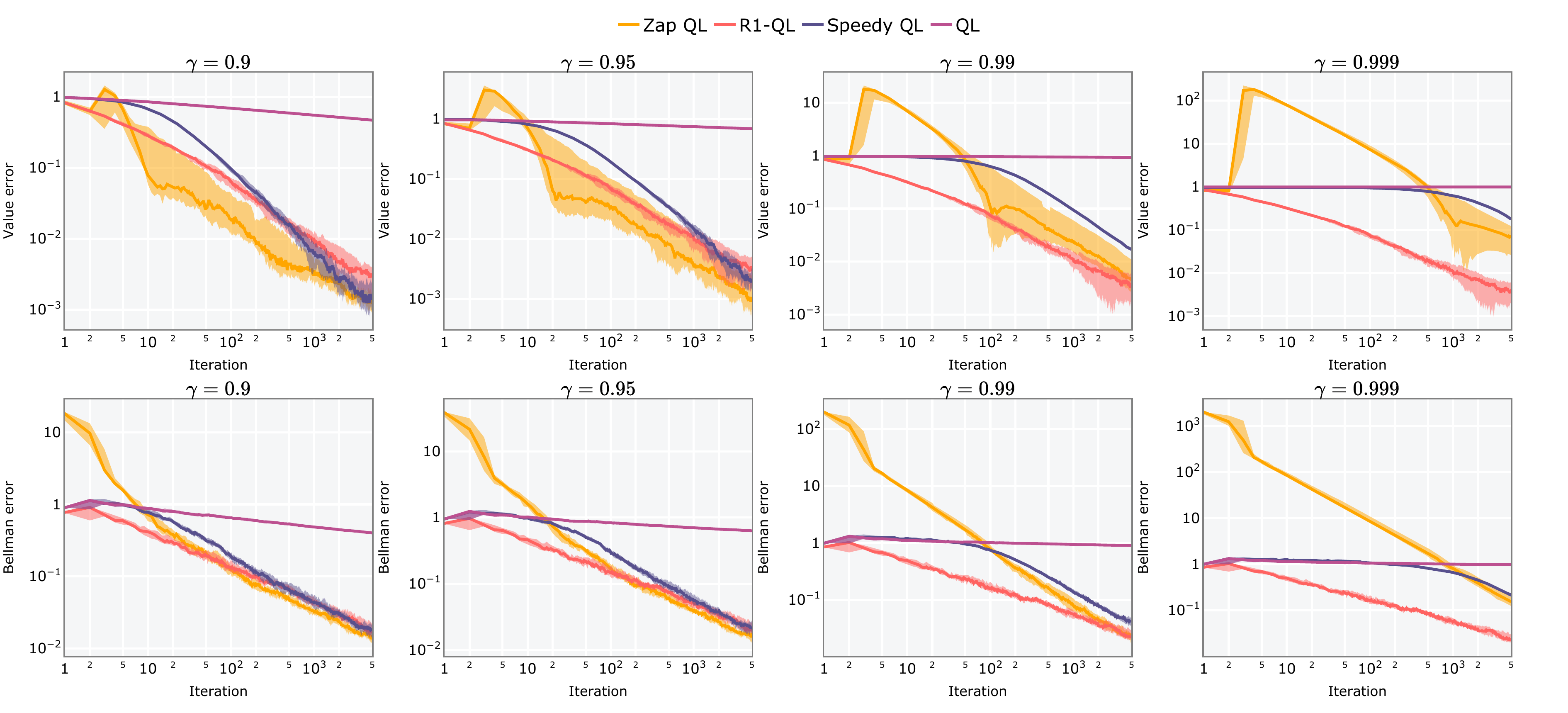}%
    \caption{Comparison of the learning algorithms in Graph MDP with various $\gamma$ values.}
    \label{fig:appx-graph-learn-curves}
\end{figure*}

We observe in Figure \ref{fig:appx-garnet-plan-curves} that Anderson VI has the highest error variance. Furthermore, the error curves for both Anderson~VI and Nesterov~VI are not monotonically decreasing, which is particularly visible for Anderson~VI. 
This is because the iterations in these algorithms are not necessarily a contraction with a guaranteed reduction in the Bellman/value error. 
Nevertheless, in our numerical experiments, both algorithms seem to be convergent. 
Figure \ref{fig:appx-graph-plan-curves} shows the error curves of planning algorithms for the Graph MDP. 
Here, the non-monotonic behavior of Anderson~VI is more apparent as the error values initially increase with an oscillating behavior. 
Similar to Garnet MDPs, R1-V1 consistently provides lower errors throughout the iterations.

Figure \ref{fig:appx-garnet-learn-curves} and \ref{fig:appx-graph-learn-curves} show the error curves for the learning algorithms in Garnet and Graph MDPs, respectively. 
In these experiments, we run each learning algorithm with 5 different seeds to marginalize the randomness in the sampling process. 
We observe that the difference between Zap QL, Speedy QL, and R1-QL is not noticeable at lower values of the discount factors (i.e., $\gamma \leq 0.95$). 
However, as the discount factor increases, particularly at $\gamma = 0.999$, the gap between the error values increases. 
At higher values of the discount factors, R1-QL consistently yields lower error values, while QL struggles to minimize the errors due to the linearly decaying step-size~$\lambda_k$. 
Furthermore, we observe that Zap QL displays higher error variance, which may be due to the inversion of the estimated ``Hessain'' $(\bm{I} - \gamma \bm{\wh{P}})$. 
In contrast, R1-QL exhibits considerably lower variance, despite implicitly performing a similar inversion. 
We argue that the lower error variance observed with R1-QL is due to the low-rank approximation of the transition probability matrix via estimation of the corresponding stationary distribution.

\subsection{Reducible MDPs}\label{app:reducible}

To illustrate R1-VI’s behavior in a reducible MDP, we replicate the comparison from Section \ref{sec:num} within the Gridworld environment of \cite{sutton2018reinforcement}, which inherently contains an absorbing state. We consider two Gridworld variants:
\begin{enumerate}
    \item \textbf{Absorbing Gridworld}, in which the absorbing state yields positive reward.
    \item \textbf{Terminal Gridworld}, in which the absorbing state grants a zero reward.
\end{enumerate}

Here, “reducibility” refers to the Markov chain induced by the optimal policy. Note, however, that even in the absence of an absorbing state, where all actions from a state lead back to itself, a non-optimal policy may still produce a reducible chain in Gridworld.

\begin{figure*}[!htb]
    \centering    \includegraphics[width=\textwidth]{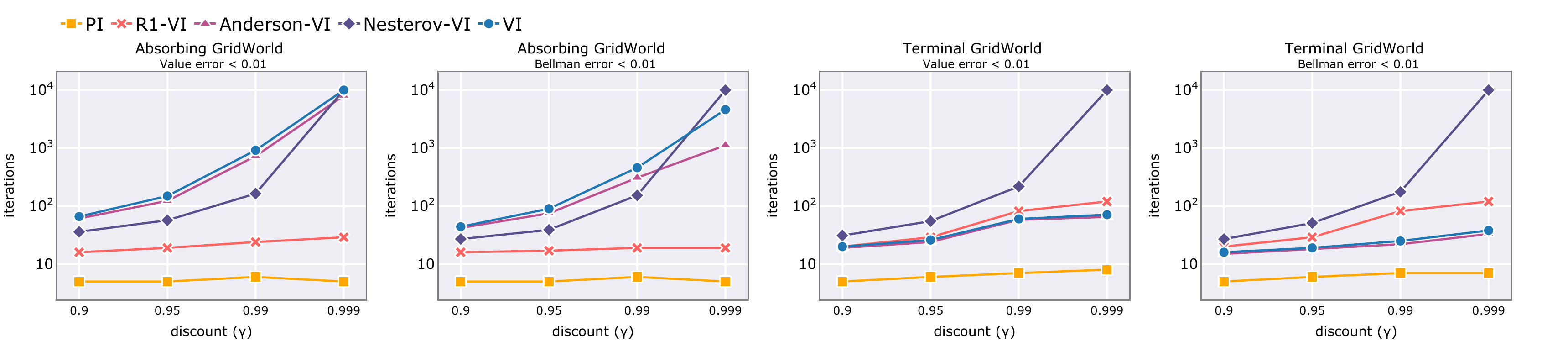}
    \caption{Comparison of planning algorithms on two reducible Gridworld MDP instances, one with a zero-reward absorbing state and one with a positive-reward absorbing state, across various $\gamma$ values.}
    \label{fig:appx-reducible-mdp}
\end{figure*}

In Absorbing-Gridworld (left side of Figure \ref{fig:appx-reducible-mdp}), R1-VI yields the lowest value and Bellman error, apart from PI, across all~$\gamma$ values. However, in Terminal Gridworld (right side of Figure \ref{fig:appx-reducible-mdp}), R1-VI performs slightly worse than the other accelerated algorithms. This is explained by the fact that under the optimal policy, the stationary distribution concentrates on the absorbing state, which provides zero reward and hence zero Bellman error when the values are initialized to zero. Consequently, the second term in the R1-VI update rule~\eqref{eq:R1-VI_first} vanishes. In contrast, in Absorbing Gridworld, the Bellman error at the absorbing state is not immediately zero, hence the second term in R1-VI contributes to improve convergence.

\subsection{Policy performance }\label{app:policy_performance}

Throughout Section \ref{sec:num}, we compared both the planning and learning algorithms, including our proposed R1VI and R1QL methods, using the value and Bellman error metrics. However, rapid convergence in value does not necessarily translate into equally rapid convergence in policy space, which is the ultimate criterion of policy optimization. In this section, we present a comparative analysis based on the policy evaluation metric in the Graph and Garnet MDPs.

\begin{figure*}[!htb]
    \centering    \includegraphics[width=\textwidth]
{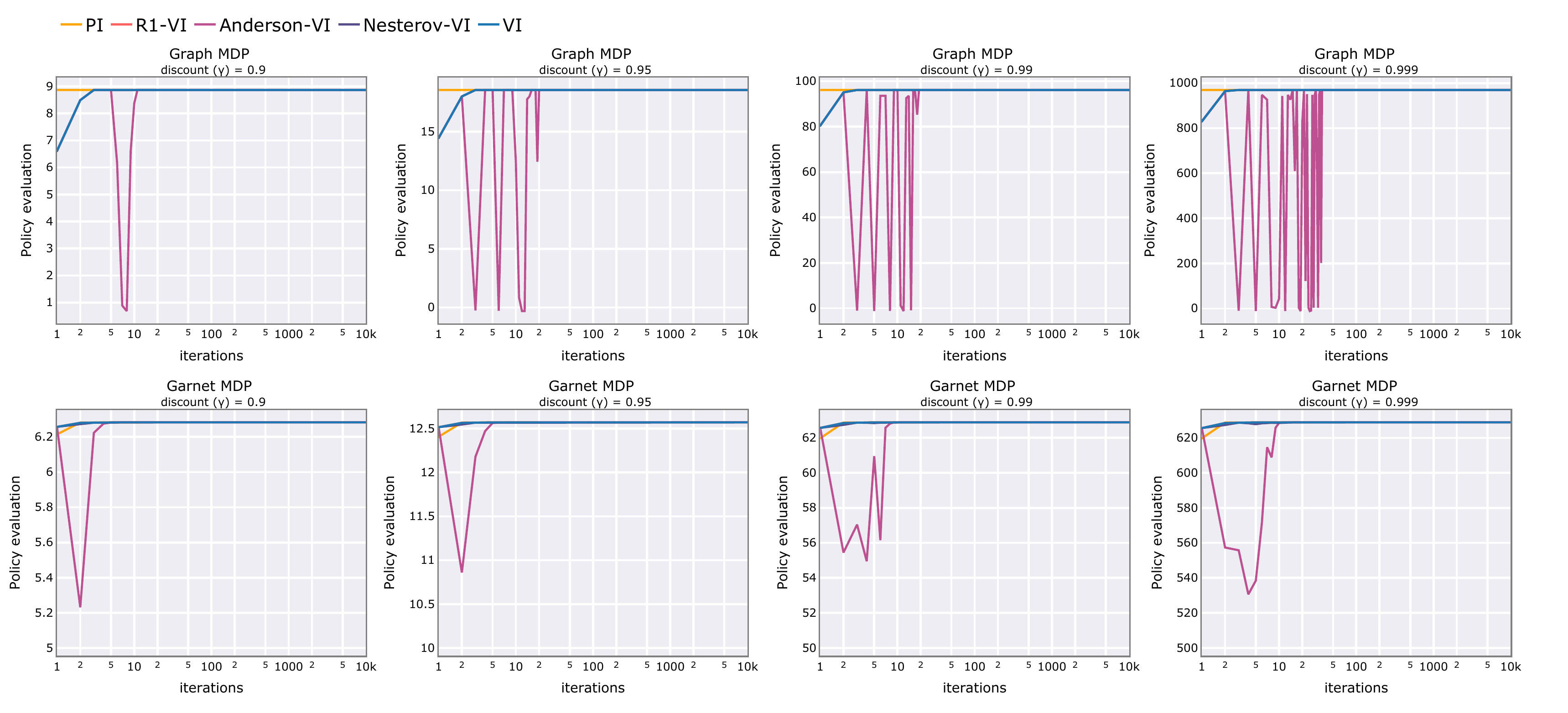}
    \caption{Comparison of the value of the greedy policy produced by various planning algorithms (R1-VI, VI, PI, Anderson-VI, and Nesterov-VI) on Garnet and Graph MDPs over a range of discount factors $\gamma$. Note that R1-VI, VI, and Nesterov-VI yield essentially overlapping results.}
    \label{fig:appx-policy_eval-plan}
\end{figure*}

\begin{figure*}[!htb]
    \centering        \includegraphics[width=\textwidth]{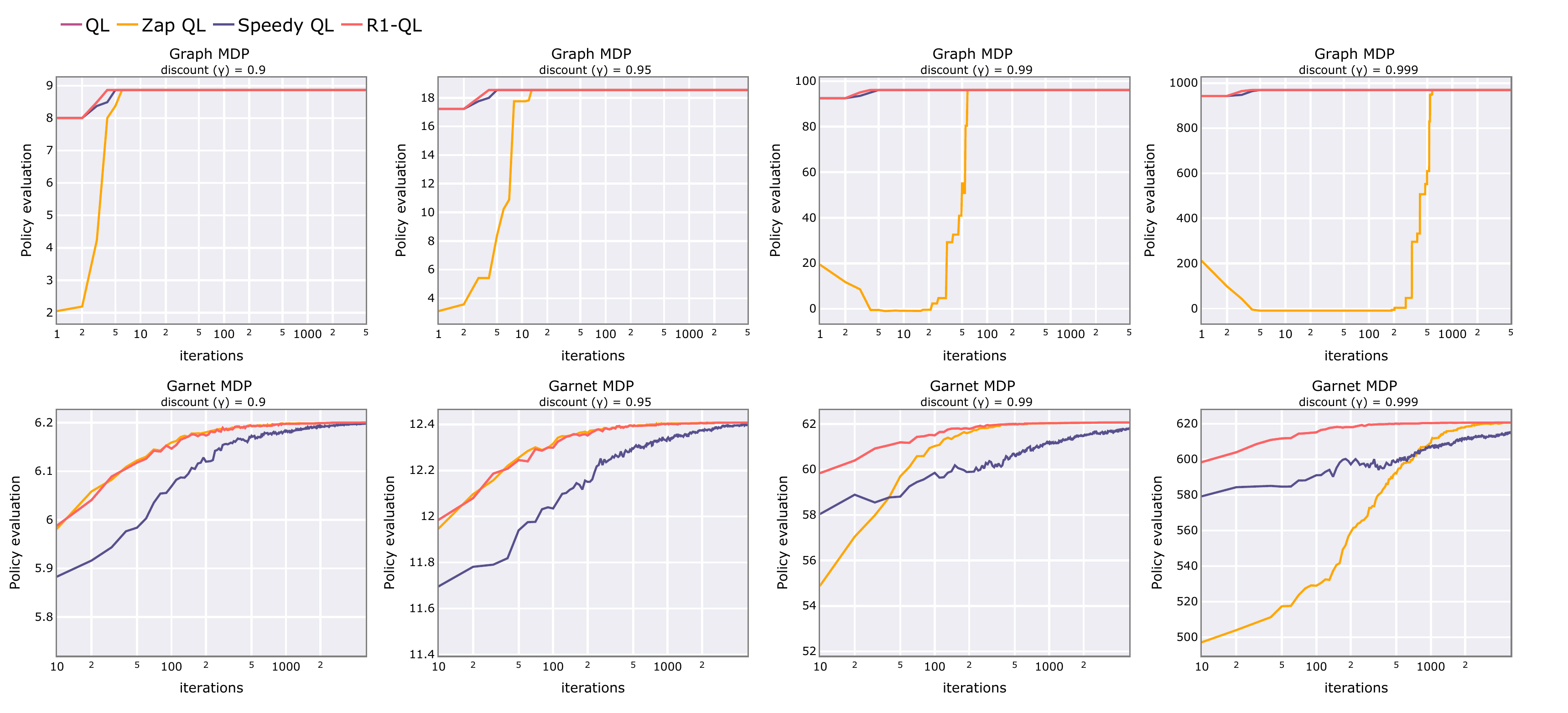}
    \caption{Comparison of the value of the greedy policy produced by the R1-QL and QL learning algorithms on Garnet and Graph MDPs as a function of the discount factor $\gamma$. Results for R1-QL and QL coincide. The y-axis shows the median policy-evaluation score of the corresponding greedy policies across 5 different seeds.}
    \label{fig:appx-policy_eval-learn}
\end{figure*}

Figure \ref{fig:appx-policy_eval-plan} shows that the planning algorithms yield exactly the same policy evaluation, except for PI and Anderson‐VI. Anderson‐VI initially struggles to find the optimal policy due to instabilities in the value space (shown in Figures \ref{fig:appx-garnet-plan-curves} and \ref{fig:appx-graph-plan-curves}), but eventually converges to the optimal policy. In both MDPs, policy convergence occurs in fewer than five steps, except for Anderson‐VI, whereas convergence in the value space requires several orders of magnitude more steps.

The convergence of policies among the learning algorithms is more varied than that of the planning algorithms. Figure \ref{fig:appx-policy_eval-learn} shows that, for the Graph MDP, all algorithms except Zap‐QL achieve almost identical performance, converging within five steps for every value of $\gamma$. In the Garnet MDP, convergence requires more steps, and improvements over iterations are slower.

In both Figures \ref{fig:appx-policy_eval-plan} and \ref{fig:appx-policy_eval-learn}, the proposed R1VI and R1QL algorithms match the policy performance of VI and QL, respectively, across all iterations. Moreover, VI in the planning setting and QL in the learning setting produce among the highest policy performance observed across both MDPs, despite exhibiting the slowest convergence in the value space.


\bibliographystyle{apalike} 
\begin{small}
\bibliography{ref}
\end{small}

\end{document}